\documentclass[a4paper,12pt,twoside]{article}
\input{Commandes.sty}

\TitleIs{Complex moments of class numbers with fundamental unit restrictions} 
\AuthorIs{Jérémy Dousselin} 

\begin{document}
\entete{We explore the distribution of class numbers $h(d)$ of indefinite binary quadratic forms, for discriminants $d$ such that the corresponding fundamental unit $\varepsilon_d$ is lower than $d^{1/2+\alpha}$, where $0<\alpha<1/2$. To do so we find an asymptotic formula for $z^{th}$-moments of such $h(d)$'s, over $d\leq x$, uniformly for a complex number $z$ in a range of the form $|z|\leq(\log x)^{1+o(1)}$, $\Re(z)\geq -1$. This is achieved by constructing a probabilistic random model for these values, which we will use to obtain estimates for the distribution function of $h(d)$ over our family. As another application, we give an asymptotic formula for the number of $d$'s such that $h(d)\leq H$ and $\varepsilon_d\leq d^{1/2+\alpha}$ where $H$ is a large real number.}
\section{Introduction}
Let $\D:=\{d\in\NN: d\equiv 0,1~\mod 4,$ and $d$ is not a square$\}$ be the set of positive discriminants. For $d\in\D$, let $h(d)$ be the number of equivalence classes of primitive binary indefinite quadratic forms $ax^2+bxy+cy^2$ of discriminant $d=b^2-4ac$. Ever since Gauss' 200 years old {\it Disquisitiones Arithmeticae}, the study of the class number $h(d)$ is still active, and many questions remain unanswered. While we know a lot about binary forms corresponding to negative discriminants, little is known about its positive counterpart. For example, Heilbronn \cite{heilbronn} proved in 1934 that $h(d)\to\infty$ when $d\to-\infty$, but it is still an open problem to know if $h(d)=1$ for infinitely many positive discriminants. The main difference between these cases is the presence of non-trivial units, in the case $d>0$, which affects the size of $h(d)$. This is explained by Dirichlet's Class Number Formula, which states that 
\[\forall d\in \D, \quad h(d)=\frac{L(1,\chi_d)\sqrt d}{\log\varepsilon_d}.\]
Here and throughout $\chi_d=\left(\frac d\cdot\right)$ is Kronecker's symbol, and $L(s,\chi_d)$ is the Dirichlet $L$-function attached to it. The fundamental unit $\varepsilon_d$ will be defined by $(t+\sqrt d u)/2$, if $(t,u)\in\NN\times\NN^*$ is the smallest solution in $t$ of the Pell equation $t^2-du^2=4$. We also define, for $A\subset \RR$, $A(x):=\{a\in A: a\leq x\}$.

Gauss already understood that the irregularities of $h(d)$ are due to huge variations of the regulator $\log \varepsilon_d$ (that may be as small as $\log d$ or as large as $\sqrt d$). For instance, instead of directly studying the mean of $h(d)$ for $d\in\D(x)$, he conjectured in 1801 that
\[\sum_{d\in \D(x)}h(d)\log\varepsilon_d\sim \frac{\pi^2}{18\zeta(3)}x^{3/2}.\]
Siegel \cite{siegel} proved this asymptotic formula in 1944, and we had to wait until 1984 for Hooley \cite[Conjecture 7]{hooley_pellian_1984} to give a conjecture -which is still an open problem to this day- concerning the mean of $h(d)$. Hooley's conjecture states that
\[\sum_{1\leq d\leq x: d\neq \square}h(d)\sim \frac{25}{12\pi^2}x(\log x)^2.\]
Here the notation $d\neq \square$ stands for the assertion "$d$ is not a perfect square."

To achieve such a {\it tour de force}, Hooley considered a family of integers with small associated fundamental unit. More precisely,  for $0\leq \alpha\leq 1$ he considered the family $\D_\alpha^{Hooley}$ of $d\in\NN$, $d\neq \square$, such that $\varepsilon_d^{Hooley}\leq d^{1/2+\alpha}$. Here $\varepsilon_d^{Hooley}$ is the fundamental unit as defined by Hooley, which slightly differs from our definition: we have $\varepsilon_d^{Hooley}=t+\sqrt du$, where $(t,u)$ is the smallest solution in $t$ of the Pell equation $t^2-du^2=1$. We will detail some consequences of this difference in \nameref{rema5} and \ref{sec:hooley_diff}. 

Hooley proved that uniformly in $0<a\leq \alpha\leq 1/2$, where $0<a<1/2$ is fixed:
\[\sum_{d\in\D_\alpha^{Hooley}(x)}1\sim \frac{4\alpha^2}{\pi^2}\sqrt x\log^2x.\]
If $0<\alpha<1/2$ is fixed, he also proved that 
\[\sum_{d\in\D_\alpha^{Hooley}(x)}h(d)\sim \frac{4(2\alpha-\log(1+2\alpha))}{\pi^2}x\log x.\]
Since the fundamental unit $\varepsilon_d^{Hooley}$ used by Hooley is defined in a different way than ours, the constant terms above will differ from ours. 

In our case, it will be more suitable to define:
\[\numberthis\label{eq:defDa}\D_\alpha:=\{d\in\D:\varepsilon_d\leq d^{1/2+\alpha}\}.\]%
Following Hooley's idea and assuming GRH, Dahl \cite{dahl_moments_2010} showed (with our definition of $\varepsilon_d$) that for $\lambda>-1$ and $0<\alpha<1/2$ fixed:
\[\sum_{d\in\D_\alpha(x)}h(d)^\lambda\sim C(\lambda,\alpha)x^{(\lambda+1)/2}\log(x)^{2-\lambda},\]
where $C(\lambda,\alpha)$ is an explicit constant. He also found a different formula in the case $\lambda=-1$. Combining ideas introduced by Lamzouri \cite{lamzouri_large_2017} and Dahl-Lamzouri \cite{dahl_distribution_2016}, we are able to improve upon this result, allowing $\lambda$ to be a complex function of $x$ (uniformly in some range $|\lambda|\leq (\log x)^{1+o(1)}$, $\Re(\lambda)\geq-1)$, and removing the need for GRH. However in order to properly state our results, we first need to construct a random model for $L(1,\chi_d)$ when $d$ varies in $\D_\alpha(x)$. 
\subsection{Settings of the problem and main statement}\label{subsec:intro}
Using techniques of Granville-Soundararajan \cite{gville} and Dahl-Lamzouri \cite{dahl_distribution_2016}, we wish to compare the distribution of $L(1,\chi_d)$ for $d\in\D_\alpha(x)$ to that of a random Euler product $L(1,\XX):=\prod_p(1-\XX(p)/p)^{-1}$, where the $\XX(p)$'s are independent random variables taking the values $0,\pm1$. The idea is to think of $\XX(p)$ as a random model for the value $\chi_d(p)$ when $d\in\D_\alpha$. To do so, we define 
\[a_p:=\PP(\XX(p)=1):=\frac{(p-1)^2}{2p(p+1)},\qquad b_p:=\PP(\XX(p)=-1):=\frac{p-1}{2p},\]
\[c_p:=\PP(\XX(p)=0):=\frac{2}{p+1}.\]
The reason behind these choices of probabilities, instead of the easier $a_p=b_p=1/2$, is the following. If $d\in\D_\alpha$, then $d=\tilde d\ell^2$, with some $\ell\in\NN^*$ and $\tilde d$ a fundamental discriminant. A fundamental discriminant being divisible by a prime $p$ with probability $|p\ZZ/p^2\ZZ-\{0\}|/|\ZZ/p^2\ZZ-\{0\}|=1/(p+1)$, the multiplicativity of $\chi_{\cdot}(p)$ implies that $c_p=1/(p+1)+(1-1/(p+1))\times 1/p=2/(p+1)$. Moreover, thanks to calculations carried out in the proof of \nameref{lemm:esp} below, we know that we should have
\[a_p-b_p=\EE(\XX(p))=-\frac{p-1}{p(p+1)}.\]
Knowing that $a_p+b_p=1-c_p=1-2/(p+1)$, we find the values of $a_p$, $b_p$ and $c_p$ above.
\begin{rema}
This probabilistic random model is to be compared with the one introduced by Granville and Soundararajan in \cite{gville}, where the probabilities that their random variables were equal to 1 and -1 were both equal to $p/(2(p+1))$. In our probabilistic random model however, this is not the case and it appears that our random variables should be equal to -1 more often than they are equal to 1. It may also be interesting to note that our probabilities are independent of $\alpha$, which  was not expected at first sight.
\end{rema}
Here is a plot comparing the probabilistic random model to the actual distribution of $\chi_d(p)$ for $d\in\D_{1/4}(10^5)$.
\begin{center}
\includegraphics[scale=0.55]{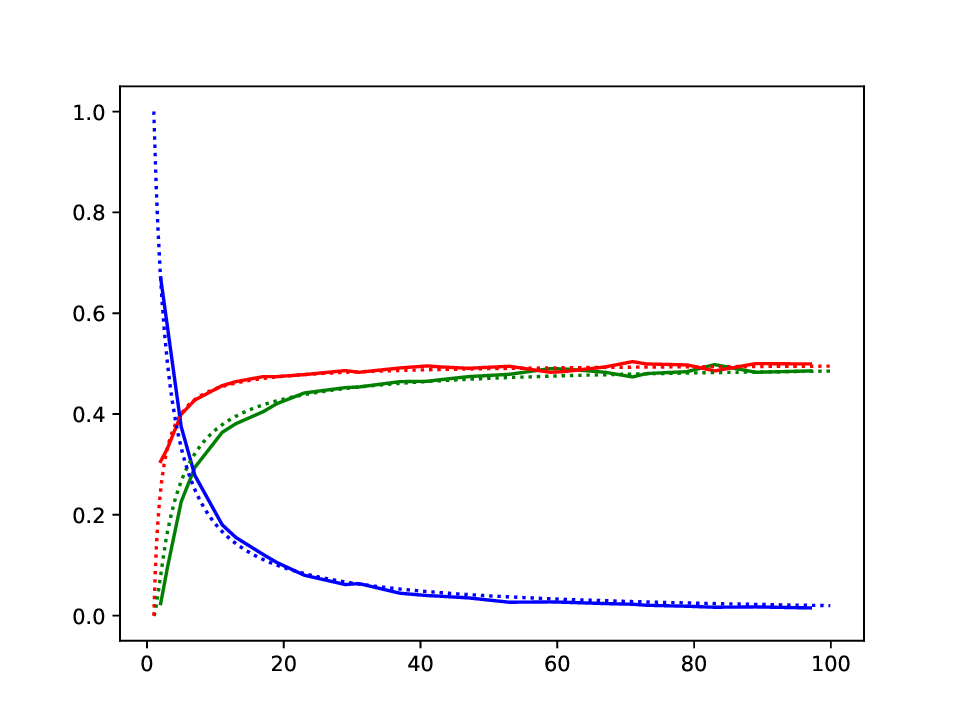}\\\it
The coloured curves are the plots of the functions that map a prime $p$ to the proportion of $d\in\D_{1/4}(10^5)$ such that $\chi_d(p)=1$ (in green)/ $\chi_d(p)=-1$ (in red)/ $\chi_d(p)=0$ (in blue). The dotted curves represent the plots of the corresponding theoretical probabilities $p\mapsto a_p,p\mapsto b_p$ and $p\mapsto c_p$.
\end{center}
\vspace{1em}
We define $\XX(n)$ by total multiplicativity, setting $\XX(1)=1$ and $\XX(p_1^{\beta_1}...p_r^{\beta_r}):=\XX(p_1)^{\beta_1}...\XX(p_r)^{\beta_r}$. 

Here and throughout, we write $\log_k$ for the $k^{th}$-iterate of the natural logarithm function. We may now state our main result. 
\begin{theo}\label{theo-hd}
Let $x$ be large and $0<\alpha<1/2$ be fixed. Uniformly for every complex number $z$ satisfying $|z|\leq \frac{(1-2\alpha)^2}{75}\log x/(\log_2 x\log_3 x)$ and $\Re(z)\geq -1$, we have that
\begin{align*}
\sum_{d\in\D_\alpha(x)}h(d)^z&=\frac{\EE(L(1,\XX)^z)}{2\zeta(2)}\int_2^x\frac{t^{\frac{z-1}2}}{\log(t)^z}\Biggl[I_{z,1}(\alpha)\log^2(t)+\phi(z)I_{z,0}(\alpha)\log(t)\Biggr]\diff t+Err
\end{align*}
Here,
\[ I_{z,1}(\alpha):=\int_{1/2}^{\alpha+1/2} u^{-z}\left(u-\frac 12\right)\diff u,\quad I_{z,0}(\alpha):=\int_{1/2}^{\alpha+1/2} u^{-z}\diff u,\] 
\[\phi(z)=\sum_{p\geq 3}\frac{\frac{2\log p}{(p+1)^2}\left(1-\left(1-\frac1p\right)^{-z}\right)}{\EE\left(\left(1-\frac{\XX(p)}{p}\right)^{-z}\right)}-\frac{\log 2}9\frac{4\times2^z+5}{\frac1{6}2^z+\frac{1}2\left(\frac23\right)^z+\frac43}-\frac{2\zeta'(2)}{\zeta(2)}+2\gamma,\]
with $\gamma$ being the usual Euler-Mascheroni constant and 
\[Err\ll \frac{2^{\Re(z)}\EE(L(1,\XX)^{\Re(z)})}{\Re(z)+1}\frac{x^{(\Re(z)+1)/2}}{(\log x)^{\Re(z)}}\left(\frac{|z|^2}{(\Re(z)+1)\log x}+\log\log^3x\left(1+\frac{|z|}{\Re(s)+2}\right)\right)\]
if $\Re(z)>-1$ and, if $\Re(z)=-1$
\[Err\ll |z|(\log^2 x)\log\log^3 x.\]
\end{theo}

Note that for some $z$ -for example any zero of $\EE((1-\XX(p)/p)^{-z})$ for $p$ a prime-, $\phi(z)$ is not properly defined, but the product $\EE(L(1,\XX)^z)\phi(z)$ is. This naturally gives a meaning to the formula above for such $z$'s.

More explicitly we have the following corollary:
\begin{coro}\label{coro_equiv}
Let $x$ be a large real number, $0<\alpha<1/2$ be fixed, and let $z$ be a complex number satisfying $|z|\leq\frac{(1-2\alpha)^2}{75}\log x/(\log_2 x\log_3 x)$. We also assume that $\Re(z)$ is large and that $\Im(z)\ll \Re(z)$. Then we have that
\[\sum_{d\in\D_\alpha(x)}h(d)^z=\frac{1}{\zeta(2)}\frac{x^{\frac{z+1}2}}{z+1}I_{z,1}(\alpha)\log(x)^{-z+2}\left(\EE(L(1,\XX)^z)+\O\left(\frac{\EE(L(1,\XX)^\sigma)}{\log_2x\log_3 x}\right)\right)+Err,\]
where $Err$ is defined in \nameref{theo-hd} above. In particular if $z$ is a large real number such that $z\leq \frac{(1-2\alpha)^2}{75}\frac{\log x}{(\log_2 x)^2}$, then
\[\sum_{d\in\D_\alpha(x)}h(d)^z=\frac{1}{\zeta(2)}\frac{\EE(L(1,\XX)^z)}{z+1}I_{z,1}(\alpha)x^{\frac{z+1}2}\log(x)^{-z+2}\left(1+\O\left(\frac1{\log_2 x}\right)\right).\]
\end{coro} 
To prove our main theorem, we have to study moments of $L$-functions associated with $\chi_d$, for $d\in\D_\alpha(x)$. 
\begin{theo}\label{prop:sum_Lz}
Let $x\geq 2$ and $0<\alpha'<1/2$ be fixed. Then uniformly for $0< \alpha\leq \alpha'$ and $|z|\leq\frac{(1-2\alpha)^2}{75}\frac{\log x}{\log_2x\log_3x}$, we have
\begin{align*}
\sum_{d\in\D_\alpha(x)}L(1,\chi_d)^z&=\frac{\sqrt x}{\zeta(2)}\EE(L(1,\XX)^z)\Biggl(\alpha^2\frac{\log^2x}{2}+\alpha\left(\phi(z)-2\alpha\right)(\log x-2)\Biggr)
\\&+\O\left(\sqrt x(\log\log x)^3\EE(L(1,\XX)^{\Re(z)})\right).
\end{align*}
Here $\phi$ was defined in \nameref{theo-hd}. 
\end{theo}
This result might be compared with the ones of \cite[Theorem 1.4]{dahl_distribution_2016} and \cite[Theorem 2]{gville}, where the authors obtained similar estimates (over different families of discriminants).
\begin{rema}\label{rema5}
Unlike the results of \cite{dahl_distribution_2016} and \cite{gville}, we encountered a major technical difficulty in the proof of \nameref{theo-hd}. In order to use \nameref{prop:sum_Lz} to prove \nameref{theo-hd}, we have to use Stieltjes integration. This forces us to obtain our \nameref{prop:sum_Lz} {\bf uniformly} in $0<\alpha<\alpha'<1/2$. When proving this result, we observe that the behaviour of the studied sum is quite different depending on the size of $\alpha$ compared to $x$. To understand why, let us recall some facts about the fundamental unit. If $\varepsilon_d=(t+u\sqrt d)/2$ where $t^2-du^2=4$, $u\geq 1$, we may deduce that $t\geq \sqrt d$ and hence $\varepsilon_d\geq (1+u)\sqrt d/2$. Therefore, when $\alpha$ is relatively small compared to $x$ ($\alpha\approx 1/\log x$ as we will show later), the condition $\varepsilon_d\leq d^{1/2+\alpha}$ implies that $u=1$. We will give more precise details in \ref{sec:hooley_diff}, but until then we can outline the following heuristic:
\begin{enumerate}
\item If $\alpha$ is very small, then $\D_\alpha(x)$ is empty;
\item If $\alpha$ is of intermediate size (so that we can only have $u=1$), then the only discriminants in $\D_\alpha(x)$ are of the form $d=t^2-4$;
\item If $\alpha$ is larger, then we may have $u>1$, and this case is the only one contributing to the main term of our theorem above.
\end{enumerate}
The second case corresponds to discriminants that are almost like those in Chowla's family, studied in \cite{dahl_distribution_2016}. We will say that the discriminants from the third case are of {\it Hooley's type}. Hooley did not encounter such difficulties with the intermediate case, for his definition of the fundamental unit $\varepsilon_d^{Hooley}=t+\sqrt du$ and $t^2-du^2=1$ implies that $\varepsilon_d^{Hooley}\geq 2\sqrt d$, and hence $\alpha$ cannot be too small. In our case, we will deal with the very difficult transition from Hooley's discriminants to Chowla's discriminants by stating our intermediary results according to the size of $\alpha$, and we will then show that the contribution of Chowla's type discriminants goes in the error term.
\end{rema}

\subsection{Study of the distribution of $h(d)$ over $d\in\D_\alpha(x)$}\label{sec:2}
\nameref{coro_equiv} might be used to study the distribution function of $h(d)$'s over our family of discriminants. We will prove the following:
\begin{theo}\label{theo:distribution}
Let $x$ be large. For some positive constant $B$, uniformly in the range $B\leq \tau \leq \log_2 x-2\log_3 x+2\log(1-2\alpha)-\log(75)$, the proportion of $d\in\D_\alpha(x)$ such that 
\[h(d)\geq 2e^\gamma \frac{\sqrt x}{\log x}\tau\]
equals 
\[\exp\left(-\frac{e^{\tau-C_0}}\tau\left(1+\O\left(\frac{1}{\sqrt\tau}\right)\right)\right),\]
where 
\[C_0:=\int_0^1\frac{\tanh(t)}t\diff t+\int_1^\infty\frac{\tanh(t)-1}t\diff t\approx 0.819.\]
\end{theo}
Here is an illustration of this fast decreasing proportion.
\begin{center}
\includegraphics[scale=0.55]{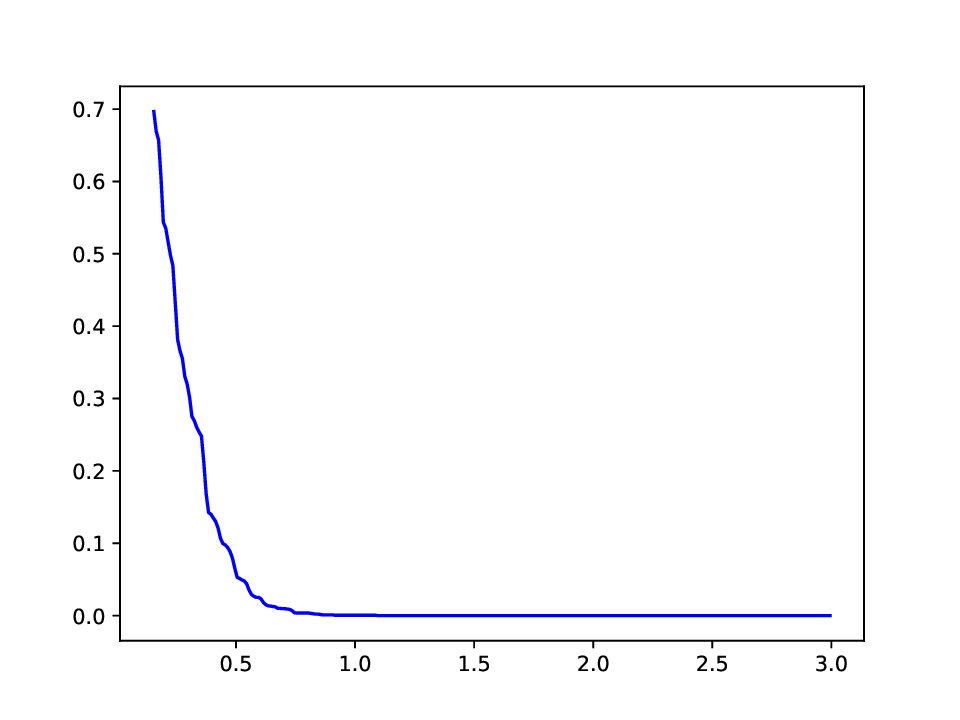}\\\it 
The graph of the function that maps $t$ to the proportion of $d\in\D_{1/10}(10^6)$ such that $h(d)\geq 2e^\gamma\frac{10^{3}}{6\log 10}t$. 
\end{center}\vspace{0.5cm}
Although their probabilistic random models differ from ours, the random variable $L(1,\XX_{Chowla})$ introduced by Dahl and Lamzouri \cite{dahl_distribution_2016} to "mimic" $L(1,\chi_d)$ when $d\leq x$ is in Chowla's family, as well as the random variable $L(1,X)$ introduced by Granville and Soundararajan \cite{gville} to approximate $L(1,\chi_d)$ when $d$ is a fundamental discriminant such that $|d|\leq x$ both satisfy a similar estimate to that of \nameref{theo:phi_esti} below (see respectively their Theorem 1.3 and Proposition 1). More precisely, the tails of the probability distribution of these values share a common "double exponential decay," as properly stated in the theorem below, in the case of our probabilistic random model.
\begin{theo}\label{theo:phi_esti}
We define $\Phi_\XX(\tau):=\PP(L(1,\XX)>e^\gamma\tau)$ and $\Psi_\XX(\tau):=\PP\left(L(1,\XX)<\frac{\zeta(2)}{e^\gamma\tau}\right)$. For large $\tau$, we have 
\[\Phi_\XX(\tau)=\exp\left(-\frac{e^{\tau-C_0}}{\tau}\left(1+\O\left(\frac1\tau\right)\right)\right).\]
The same estimate holds for $\Psi_\XX$. Moreover, if $0\leq \lambda\leq e^{-\tau}$, then we have
\[\Phi_\XX(e^{-\lambda}\tau)=\Phi_\XX(\tau)(1+\O(\lambda e^{\tau})),\text{ and }\Psi_\XX(e^{-\lambda}\tau)=\Psi_\XX(\tau)(1+\O(\lambda e^{\tau})).\]
\end{theo}
The proof is the same to that of \cite[Theorem 1.3]{dahl_distribution_2016}, since our probabilities satisfy $a_p,b_p=1/2+\O(1/p)$ and $a_p-b_p, c_p\ll 1/p$.

\subsection{Counting the number of Hooley's quadratic fields with bounded class numbers}
Another application of our main theorem is given in \ref{sec:number_of}, where we aim to find an asymptotic formula in $H$ for $\sum_{h\leq H}\F_\alpha(h)$, the number of $d$'s in $\D_\alpha$ such that $h(d)\leq H$. Here $\F_\alpha(h)$ is defined as the number of $d\in\D_\alpha$ such that $h(d)=h$. The method we use was designed by Soundararajan \cite{soundararajan}, and he used it to find an asymptotic formula for the average value of the number of imaginary quadratic fields with class number $h$. Heuristically, we know that $L(1,\chi_d)$ is constant on average, thanks to \nameref{prop:sum_Lz}. Thus the main contribution to the average we want to estimate comes from $d\ll H^2(\log H)^2$, by the Class Number Formula. Since there are around $H\log^3 H$ such discriminants by \nameref{prop:sum-chid}, we may expect the average of $\F_\alpha$ to be around $H\log^3H$. This is indeed the case.
\begin{theo}\label{theo:number_F}
We fix $0<\alpha<1/2$. Let $H$ be a positive integer and $\varepsilon>0$ be small. Then, as $H\to\infty$, we have that
\[\sum_{h\leq H}\F_\alpha(h)=2\frac{\alpha^2(4\alpha+3)}{3}\prod_p\left(1-\frac{2p-1}{p^4}\right)H\log^3 H\left\{1+\O_\varepsilon\left(\frac{1}{(\log H)^{1/3-\varepsilon}}\right)\right\}.\]
\end{theo}
\begin{rema}
If we had uniformity in the parameter $\alpha$ and choose $\alpha\asymp 1/\log H$, the sum above should roughly be the number of Chowla's discriminants with class number $\leq H$. Our Theorem would show that this quantity should be $\asymp H\log H$, which is precisely the result obtained by Dahl and Lamzouri \cite{dahl_distribution_2016}: if $\F_{\mathrm{ch}}(h)$ is the number of discriminants $d$, squarefree, of the form $d=4m^2+1$ for some $m\geq 1$, and with class number $h$, then
\[\sum_{h\leq H}\F_{\mathrm{ch}}(h)=\frac{1}{2G}H\log H\left(1+\O\left(\frac{(\log_2 H)^2\log_3 H}{\log H}\right)\right),\]
where $G=L(2,\chi_{-4})$ is Catalan's constant, and $\chi_{-4}$ is the non-principal character modulo $4$.
\end{rema}

One can also compare this with the same result in the context of imaginary quadratic fields, proved by Soundararajan \cite{soundararajan}: if $\F(h)$ is the number of imaginary quadratic fields with class number $h$, then, as $H\to\infty$,
\[\sum_{h\leq H}\F(h)=\frac{3\zeta(2)}{\zeta(3)}H^2\left(1+\O\left((\log H)^{-1/2+\varepsilon}\right)\right).\]

Due to weaker error terms than theirs in our formula for complex moments of class numbers, our proof will be slightly more technical than theirs.\vspace{0.3cm}

\noindent\underline{\bf Acknowledgement}: I would like to thank my PhD advisor Youness Lamzouri for his time, guidance, and several useful suggestions through the preparation of this paper.

\section{On the differences of context between Hooley's paper and ours}\label{sec:hooley_diff}
In his paper, Hooley \cite{hooley_pellian_1984} considers the number $h(d)$ of properly primitive classes of indefinite forms $ax^2+2bxy+cy^2$. However as mentioned by Davenport at the beginning of the section 6 of \cite{davenport}, Lagrange and now most modern writers would instead use the notation $ax^2+bxy+cy^2$ for a quadratic form (see \cite{sarnak} for instance). 

Hooley then defines the discriminant of his quadratic form to be $d=b^2-ac$, while ours is defined by $d=b^2-4ac$. In other words, while we define the set of positive discriminants by
\[\D=\{d\in\NN : d\equiv 0\mod 4 \text{ or }d\equiv 1\mod 4; \quad d\neq \square\},\]
Hooley only studies $d/4$ when $d\in \D$ and $d\equiv 0\mod 4$. By doing so Hooley avoids the technical difficulties arising from the case $d\equiv 1\mod 4$, described in \nameref{rema5} above. 

We shall now start to give some lemmas, and we will highlight how Hooley's definition would change some of them. First, we start by parametrizing our family $\D_\alpha(x)$ defined in \eqref{eq:defDa}, $x>0$ large and $0<\alpha<1/2$. For $t\geq 0$, $u\geq 1$, we define 
\[d(t,u):=\frac{t^2-4}{u^2}.\]
As shown by Hooley, we have the following lemma:
\begin{lemm}\label{lemm}
Let $x$ be large, $0<\alpha<1/2$, and put $X_\alpha:=x^{\alpha}-x^{-1-\alpha}$. Let $Y_1:=Y_1(u,\alpha)$ be the unique positive solution to $u=Y_1^\alpha-Y_1^{-1-\alpha}$, $Y_2=(Y_1u^2+4)^{1/2}$ and $Y_3=(xu^2+4)^{1/2}$. Then each $d\in\D_\alpha(x)$ may be uniquely written as a $d(t,u)$, where 
\[1\leq u\leq X_\alpha,\qquad Y_2\leq  t\leq Y_3, \qquad d(t,u)\in\D.\]
\end{lemm}
\begin{proo}
Since the only unit lower that $d$ is the fundamental unit, we can write that $d\in\D_\alpha(x)$ if, and only if, 
\[\numberthis\label{eq:conditions} t^2-du^2=4, \quad\varepsilon_d=\frac{t+\sqrt d u}2\leq d^{1/2+\alpha}, \quad u\geq 1, \quad d\leq x.\]
Observe that the first condition ensures that $d\neq \square$. Furthermore, knowing that $u=\frac{\varepsilon_d-\varepsilon_d^{-1}}{\sqrt d}$, we may replace the second and fourth conditions above by 
\[u\leq d^{\alpha}-d^{-1-\alpha}\leq x^{\alpha}-x^{-1-\alpha}=X_\alpha,\]
since $d^{\alpha}-d^{-1-\alpha}$ is an increasing function of $d$. Now this inequality between $u$ and $d$ implies that $d\geq Y_1$. Furthermore the first condition of \eqref{eq:conditions} may be replaced by the following inequality for each $u$:
\[Y_2=(Y_1u^2+4)^{1/2}\leq t\leq (xu^2+4)^{1/2}=Y_3.\]
This concludes the proof.
\end{proo}
With his definition, Hooley gets a slightly different parametrization, and in particular $X_\alpha^{Hooley}=\frac12 X_\alpha$. It is now time to give a more quantitative version of \nameref{rema5} above, where the technical difficulty mentioned is related to the size of $X_\alpha$. We shall define $\alpha_0$ to be such that $X_{\alpha_0}=1$, and $\alpha_1$ to be such that $X_{\alpha_1}=2$. Observe that $\D_\alpha(x)$ is empty if $\alpha<\alpha_0$, since we would have $X_\alpha<1$. Also note that 
\[\numberthis\label{eq:def_alpha}\alpha_1\asymp 1/\log x.\]
By definition, the family studied by Hooley is empty as soon as $X_\alpha<2$, which is the case whenever $\alpha<\alpha_1$.

For several reasons detailed later, we will have to study $\sum_{d\in\D_\alpha(x)}\chi_d(m)$, for a fixed $m\geq 1$. We may use this \nameref{lemm} to write that
\begin{align*}
\sum_{d\in\D_\alpha(x)}\chi_d(m)&=\sum_{u=1}^{X_\alpha}\sum_{\substack{Y_2\leq t\leq Y_3\\ d(t,u)\in\D}}\left(\frac{d(t,u)}m\right)=\sum_{u=1}^{X_\alpha}\sum_{2<a\leq 4u^2+2}\sum_{\substack{Y_2\leq t\leq Y_3\\ d(t,u)\in\D\\ t\equiv a~\mod{4u^2}}}\left(\frac{d(t,u)}m\right).
\end{align*}
We define
\[C_{m,a,u}:=\sum_{0<\ell\leq m}\left(\frac{16u^2\ell^2+8a\ell+d(a,u)}m\right)\]
and $\rho(u)$ as the number of $2<a\leq 4u^2+2$ such that $d(a,u)\in\D$. If $t\equiv a~\mod{4u^2}$ then $d(t,u)=16u^2\ell^2+8a\ell+d(a,u)$, for some $\ell\in\ZZ$, and hence $d(a,u)$ and $d(t,u)$ are equal modulo $8$\footnote{The fact that we are computing modulo $8$ is actually a natural thing to do due to the definition of Kronecker's character.}. Therefore if $t\equiv a~\mod{4u^2}$ and $t\geq Y_2$, we have that $d(t,u)\in\D$ if, and only if, $d(a,u)\in\D$. Thus
\begin{align*}
\sum_{d\in\D_\alpha(x)}\chi_d(m)&=\sum_{u=1}^{X_\alpha}\sum_{\substack{2<a\leq 4u^2+2\\d(a,u)\in\D}}\sum_{\substack{(Y_2-a)/4u^2\leq n\leq (Y_3-a)/4u^2}}\left(\frac{d(a+4u^2n,u)}m\right)
\\&=\sum_{u=1}^{X_\alpha}\sum_{\substack{2<a\leq 4u^2+2\\d(a,u)\in\D}}\sum_{0<\ell\leq m}\left(\frac{d(a+4u^2\ell,u)}m\right)\sum_{\substack{(Y_2-a)/4u^2\leq n\leq (Y_3-a)/4u^2\\ n\equiv \ell~\mod m}}1
\\&=\sum_{u=1}^{X_\alpha}\sum_{\substack{2<a\leq 4u^2+2\\d(a,u)\in\D}}C_{m,a,u}\left(\frac{Y_3-Y_2}{4u^2m}+\O(1)\right)
\\&\numberthis\label{main_term}=\sum_{u=1}^{X_\alpha}\frac{\sqrt x-\sqrt{Y_1}}{4u}\sum_{\substack{2<a\leq 4u^2+2\\d(a,u)\in\D}}\frac{C_{m,a,u}}m+\O\left(m\sum_{u=1}^{X_\alpha}\rho(u)\right),
\end{align*}
since $|C_{m,a,u}|\leq m$. We shall require the following lemma, due to Lamzouri.
\begin{lemm}[Lemma 3.3 of \cite{lamzouri_large_2017}]\label{esti_ca}
Let $a,u$ be positive integers such that $d(a,u)\in\D$. Let $m$ be a positive integer and write $m=2^{e_1}p_2^{e_2}\cdots p_r^{e_r}$, where $e_1\geq 0$ and $e_j\geq 1$ for $2\leq j\leq r$. Also let $m_0$ be the squarefree part of $m/2^{e_1}$. Then if $(m_0,u)>1$ we have $C_{m,a,u}=0$. Moreover, if $(u,m_0)=1$ then we have 
\[\frac{C_{m,a,u}}m=b_{a,u}(m)\frac{(-1)^{\omega(m_0)}}{m_0}\prod_{\substack{2\leq j\leq r\\ e_j \text{ is even}}}\left(1-\frac 2{p_j}\right)\prod_{\substack{p_j|u\\ e_j\text{ is even}}}\left(1+\frac{1}{p_j-2}\right)\]
where $\omega(m_0)$ is the number of prime factors of $m_0$, and $b_{a,u}(m)=1$ if $m$ is odd, and if $m$ is even then
\[b_{a,u}(m)=\left\{\begin{array}{cl}
0 & \text{if }d(a,u)\equiv 0\mod 4,
\\ 1 & \text{if }d(a,u)\equiv 1\mod 8,
\\ (-1)^{e_1}& \text{if }d(a,u)\equiv 5\mod 8.
\end{array}\right.\]
\end{lemm}
At this point, there is would be no significant difference if we were to prove the same lemma in the context of Hooley's work (we would have to change $\equiv 0\mod 4$ by $\equiv 0\mod 2$, $\equiv 1\mod 8$ by $\equiv \pm 1\mod 8$, and $\equiv 5\mod 8$ by $\equiv \pm 5\mod 8$ in the expression above). A fundamental difference between Hooley's definition and ours lies in the following lemma.
\begin{lemm}\label{lemm:def_bm}
We keep the notations of the lemma above and we put $B_m(u):=\sum_{\substack{3\leq a\leq 4u^2+2\\ d(a,u)\in\D}}b_{a,u}(m)$, and let $\eta(u)$ be the number of odd prime factors of $u$. Let $r_1$ be the $2$-adic valuation of $u$. Then if $e_1=0$
\[B_m(u)=\left\{\begin{array}{cl}
4\times 2^{\eta(u)}&\text{ if }r_1=0\text{ or }1,
\\8\times 2^{\eta(u)}&\text{ if } r_1=2,
\\16\times 2^{\eta(u)}&\text{ if }r_1\geq 3.
\end{array}\right.\]
If $e_1\geq 1$, we have
\[B_m(u)=\left\{\begin{array}{cl}
2(-1)^{e_1}2^{\eta(u)}& \text{ if }r_1=0,
\\4(1+(-1)^{e_1})2^{\eta(u)}&\text{ if } r_1\geq 3,
\\0& \text{ otherwise}.
\end{array}\right.\]
\end{lemm}
\begin{proo}
It is an immediate consequence of Lemma 3.4 of \cite{lamzouri_large_2017}, as detailed in the proof of Theorem 3.1 of the said paper.
\end{proo}
The mentioned Lemma 3.4 of \cite{lamzouri_large_2017} is an estimate of incongruent solutions of some equations modulo 4 and 8. These equations are, for the sake of an example, of the form
\[d(a,u)\equiv 0 \mod 4,\]
to be solved in $a$. By the Chinese remainder Theorem, we have to solve $a^2-4\equiv 0 \mod{2^{2+2r_1}}$, and in Hooley's case we would rather have to solve $a^2-1\equiv 0 \mod{2^{2+2r_1}}$, which does not have the same number of incongruent solutions if $r_1\geq 2$. If we were to define $B^{Hooley}_m(u)$ in a similar way to that of $B_m(u)$, then we would have, if $e_1=0$:
\[B_m^{Hooley}(u)=\left\{\begin{array}{cl}
4\times 2^{\eta(u)}&\text{ if }r_1=0,
\\8\times 2^{\eta(u)}&\text{ if } r_1=1,
\\16\times 2^{\eta(u)}&\text{ if }r_1\geq 2.
\end{array}\right.\]
If $e_1\geq 1$, we would have
\[B_m^{Hooley}(u)=\left\{\begin{array}{cl}
(1+(-1)^{e_1})2^{\eta(u)} & \text{ if }r_1=0,\\
0& \text{ if }r_1=1,\\
4(1+(-1)^{e_1})2^{\eta(u)}&\text{ if } r_1\geq 2.
\end{array}\right.\]

We keep the notations of the two lemmas above. Again, we leave the technical details for later, but to estimate the main term of \eqref{main_term} we will rely on Perron's formula, and hence we shall compute the Dirichlet series (for $\Re(s)>1$):
\[\sum_{\substack{u\geq 1\\ (u,m_0)=1}}\frac{B_m(u)}{u^s}\prod_{\substack{p_j|u\\e_j\text{ is even}}}\left(1+\frac{1}{p_j-2}\right).\]
In order to get a multiplicative function of $u$ inside the sum, we should normalize $B_m(u)$ as follows.
\begin{lemm}\label{lemm:function_perron}
We keep the notations of \nameref{esti_ca}. We put $a(m)=4$ if $m$ is odd and $a(m)=2(-1)^{e_1}$ otherwise. Then for $\Re(s)>1$:
\[\Sigma(m,s):=\sum_{\substack{u\geq 1\\ (u,m_0)=1}}\frac{B_m(u)/a(m)}{u^s}\prod_{\substack{p_j|u\\e_j\text{ is even}}}\left(1+\frac{1}{p_j-2}\right)=G_m(s)\zeta(s)^2\]
where 
\[\numberthis\label{eq:def-G}G_m(s):=\kappa(m,s)\frac{\prod_{\substack{2\leq j\leq r\\e_j~even}}\left(1+\frac{2(p_j-1)}{(p_j-2)(p_j^s-1)}\right)}{\prod_{p|2m}\left(1+\frac{2}{p^s-1}\right)}\frac{1}{\zeta(2s)}\]
and 
\[\numberthis\label{eq:def-kappa}\kappa(m,s):=\left\{\begin{array}{cl}
1+\frac{2(1+(-1)^{e_1})}{4^s(2^s-1)} & \text{if }e_1\geq 1,
\\1+\frac1{2^s}+\frac{2}{4^s}+\frac{4}{4^s(2^s-1)} & \text{otherwise,}
\end{array}\right.\]
\end{lemm}
\begin{proo}
The proof of Theorem 3.1 of \cite{lamzouri_large_2017} shows that 
\[\Sigma(m,s)=\kappa(m,s)\prod_{\substack{2\leq j\leq r\\e_j~even}}\left(1+\frac{2(p_j-1)}{(p_j-2)(p_j^s-1)}\right)\prod_{p\nmid 2m}\left(1+\frac{2}{p^s-1}\right).\]
The lemma follows by noticing that for $p\geq 2$ an integer,
\[1+\frac2{p^s-1}=\frac{1-\frac{1}{p^{2s}}}{\left(1-\frac1{p^s}\right)^2}.\]
\end{proo}

Now, because of \nameref{esti_ca} and \eqref{main_term} we see that it should be interesting to define the function
\[\numberthis\label{defH} H_m(s):=\frac{(-1)^{\omega(m_0)}}{4m_0}\prod_{\substack{2\leq j\leq r\\ e_j~even}}\left(1-\frac2{p_j}\right)a(m)G_m(s).\]
In Hooley's case, this function $H_m$ would be almost the same, but we would have 
\[H_m^{Hooley}(s)=\frac{(-1)^{\omega(m_0)}}{4m_0}\prod_{\substack{2\leq j\leq r\\ e_j~even}}\left(1-\frac2{p_j}\right)\kappa^{Hooley}(m,s)\frac{\prod_{\substack{2\leq j\leq r\\e_j~even}}\left(1+\frac{2(p_j-1)}{(p_j-2)(p_j^s-1)}\right)}{\prod_{p|2m}\left(1+\frac{2}{p^s-1}\right)}\frac{1}{\zeta(2s)},\]
where 
\[\kappa^{Hooley}(m,s)=\left\{\begin{array}{cl}
4\frac{2^s}{2^s-1}\left(1+\frac 1{2^s}+\frac{2}{4^s}\right) & \text{if } e_1=0,
\\ (1+(-1)^{e_1})\left(1+\frac{4}{2^s(2^s-1)}\right) & \text{if } e_1\geq 1.
\end{array}\right.\]
In our proof of \nameref{prop:sum-chid}, we will establish that for a fixed $0<\alpha<1/2$ and fixed $m\geq 1$,
\[\sum_{d\in\D_\alpha(x)}\chi_d(m)\sim \frac{H_m(1)}2 \alpha^2 \sqrt x \log^2x,\]
and hence our probabilistic random model must satisfy
\[\EE(\XX(m))=\lim_{x\to \infty}\frac{\sum_{d\in\D_\alpha(x)}\chi_d(m)}{\sum_{d\in\D_\alpha(x)}1}=\frac{H_m(1)}{H_1(1)}=H_m(1)\zeta(2).\] 
One would find a similar equivalent with Hooley's definition, as one can observe when we see the expression for $H_1^{Hooley}(s)$ which appears between (16) and (17) of \cite{hooley_pellian_1984}. Therefore, we can compute what would be a good probabilistic random model for his definition. If we were to define random independent variables $(\XX^{Hooley}(p))$ taking values in $\{\pm 1, 0\}$, to "mimic" the behaviour of $\chi_d(p)$ for $d\in\D_\alpha^{Hooley}(x)$\footnote{Recall that $\D_\alpha^{Hooley}(x)=\{d\neq \square,\, d\leq x, \,\varepsilon_d^{Hooley}\leq d^{1/2+\alpha}\}$.}, then we a good way to define them should be
\[\PP(\XX^{Hooley}(p)=1)=a_p,\qquad \PP(\XX^{Hooley}(p)=-1)=b_p\qquad \PP(\XX^{Hooley}(p)=0)=c_p,\]
for any odd prime $p$, and 
\[\PP(\XX^{Hooley}(2)=1)=3/16, \qquad \PP(\XX^{Hooley}(2)=-1)=3/16,\qquad\PP(\XX^{Hooley}(2)=0)=5/8.\]
Similarly to our $\XX$'s, we would then put $\XX^{Hooley}(1)=1$ and for $n\geq 2$, $\XX^{Hooley}(n):=\prod_{p^e||n}\XX^{Hooley}(p)^e$. 

This is very coherent with the following heuristic: Hooley studies integers $d$ such that $4d\in \D$. Therefore, the probability that $d$ is divided by an odd prime $p$ is be the same as the probability of $p$ dividing $4d$. However, in the case $p=2$, we are looking at discriminants divisible by $8$, knowing that they are already divisible by $4$. Using Bayes' Theorem, we should have 
\[c_2^{Hooley}=\frac{\PP( 8|d; d\in\D)}{\PP(4|d;d\in\D)}.\]
Recall that a discriminant $d$ can be uniquely written as $d=\tilde d\ell^2$ with $\tilde d$ a fundamental discriminant and $\ell\geq 1$ and integer. Also recall that the fundamental discriminant $\tilde d$ must satisfy $\tilde d\equiv 1\mod 4$ or $\tilde d \equiv 8,12\mod{16}$. Therefore to have $8|\tilde d \ell^2$, we have to be in one of the three cases:
\begin{enumerate}
\item $8|\tilde{d}$, which happens for $1$ of the $6$ available classes modulo $16$ for $\tilde d$;
\item $4|\tilde d$ and $8\nmid \tilde d$, which happens for only $1$ class modulo $16$, and in this case $8|\tilde{d}\ell^2\iff 2|\ell$ (the latter happens with probability $1/2$);
\item $4\nmid \tilde d$, which happens for $4$ classes modulo $16$, and in this case $8|\tilde d\ell^2\iff 4|\ell$ (the latter happens with probability $1/4$).
\end{enumerate}
By doing similar computations for the case $4|d$, we find that the probability of a Hooley discriminant being divisible by $2$ is given by
\[\frac{1/6+1/6\times 1/2+4/6\times 1/4}{2/6+4/6\times 1/2}=5/8=c_2^{Hooley}\]
as expected.

We can go even further, but we will avoid technicalities and only outline a very informal proof. In $\textsection6$ of \cite{hooley_pellian_1984}, Hooley writes the class number formula as $h(d)=2\sqrt d L_d(1)/\log(\varepsilon^{Hooley}_d)$, where $L_d(1)$ is a Dirichlet's $L$-function satisfying $L_d(1)=\sum_{m \text{ odd}}\frac{\left(\frac{d}m\right)}m$. Therefore, by removing the contribution of the powers of $2$ in the Euler product for $L_d(1)$, we should have 
\[\sum_{d\in\D_\alpha^{Hooley}(x)}L_d(1)\sim \frac{\EE(L(1,\XX^{Hooley}))}{\EE\left(\left(1-\frac{\XX^{Hooley}(2)}2\right)^{-1}\right)}\sum_{d\in\D_\alpha^{Hooley}(x)}1\sim \frac{4}{\pi^2}\alpha^2 \sqrt x\log^2 x.\]
Then, a partial summation leads to 
\[\sum_{d\in\D_\alpha^{Hooley}(x)}\frac{2\sqrt d L_d(1)}{\log d}\sim 2\int_2^x \frac{\sqrt t}{\log t}\times\frac{4\alpha^2}{\pi^2}\frac{\log^2(t)}{2\sqrt t}\diff t\sim \frac{4\alpha^2}{\pi^2}x\log x .\]
Finally, we would get by the class number formula and Stieltjes integration:
\begin{align*}
\sum_{d\in\D_\alpha^{Hooley}(x)}h(d)&=\sum_{d\in\D_\alpha^{Hooley}(x)}\frac{2\sqrt d L_d(1)}{\log d}\frac{\log d}{\log\varepsilon_d^{Hooley}}\sim\int_0^\alpha (1/2+u)^{-1}\times \frac{8u}{\pi^2}x\log x\diff u
\\&\sim \frac{4[2\alpha-\log(2\alpha+1)]}{\pi^2}x\log x,
\end{align*}
as Hooley found.

\section{Lemmata and definitions}
\begin{lemm}\label{lemm:size_Y} Uniformly in $0<\alpha<1/2$, we have for all $u\geq 2$
\[\sqrt{Y_1}=u^{1/(2\alpha)}+\O\left(\frac1{ u^{2}}\right).\]
\end{lemm}
\begin{proo}
We have that
\[Y_1^\alpha=u+Y_1^{-1-\alpha}=u+(u+Y_1^{-1-\alpha})^{(-1-\alpha)/\alpha}=u+\O(u^{(-1-\alpha)/\alpha}).\]
Since $u^{-1/(2\alpha)}/\alpha$ is uniformly bounded for $0<\alpha<1/2$, we deduce that
\begin{align*}
\sqrt{Y_1}&=u^{1/(2\alpha)}\exp\left(\frac1{2\alpha}\log\left(1+\O(u^{-1/\alpha-2})\right)\right)=u^{1/(2\alpha)}\left(1+\O\left(\frac{1}{\alpha}u^{-1/\alpha-2}\right)\right)
\\&=u^{1/(2\alpha)}+\O(u^{-1/(2\alpha)-2}/\alpha)=u^{1/(2\alpha)}+\O(u^{-2}).
\end{align*}
\end{proo}

We shall now introduce several notations. First, we define another sequence of independent random variables, $(\XX_s(p))_p$, for $s>1/2$ and where $p$ runs over the set of primes. They are defined by the following probabilities, for $p$ an odd prime number,
\[a_p(s):=\PP(\XX_s(p)=1):=\frac{p-3}{2p}+\frac{2}{p^s+1}\frac1p,\qquad b_p(s):=\PP(\XX_s(p)=-1):=\frac{p-1}{2p},\]
\[\numberthis\label{eq:def_probas} c_p(s):=\PP(\XX_s(p)=0):=1-\frac{1}p\left(p-2+\frac2{p^s-1}\right).\]
In the case $p=2$, we have 
\[a_2(s):=\PP(\XX_s(2)=1):=\frac{1}{8^s+2^s+2},\qquad b_2(s):=\PP(\XX_s(2)=-1):=\frac12\frac{8^s-4^s+2}{8^s+2^s+2},\]
\[c_2(s):=\PP(\XX_s(2)=0):=\frac{2^{s-1}(4^s+2^s+2)}{8^s+2^s+2}.\]
We then extend the definition of these variables to integers by putting $\XX_s(n)=\prod_{p^e||n}\XX_s(p)^e$ and $\XX_s(1)=1$. Note that by definition $\XX(m)$ has the same law as $\XX_1(m)$. In addition, we define $(\XX_\infty(m))_m$ to be the random variables defined by 
\[a_p(\infty):=\PP(\XX_\infty(p)=1):=\frac{p-3}{2p},\qquad b_p(\infty):=\PP(\XX_\infty(p)=-1):=\frac{p-1}{2p},\]
\[c_p(\infty):=\PP(\XX_\infty(p)=0):=\frac2p.\]
for $p$ an odd prime and
\[a_2(\infty):=\PP(\XX_\infty(2)=1):=0,\qquad b_2(\infty):=\PP(\XX_\infty(2)=-1):=\frac12,\]
\[c_2(\infty):=\PP(\XX_\infty(2)=0):=\frac12.\]
We also put $\XX_\infty(n):=\prod_{p^e||n}\XX_\infty(p)^e$, and $\XX_\infty(1)=1$.

To isolate the multiplicative part of $H_m(s)$ defined in \eqref{defH}, we write it in the form 
\begin{align*}
 H_m(s)&=\frac{8^s+2^s+2}{8^s+4^s}\frac{1}{\zeta(2s)}\Biggl\{\frac{(-1)^{\omega(m_0)}}{m_0}\frac{\prod_{\substack{2\leq j\leq r\\ e_j~even}}\left(1-\frac2{p_j}\right)\left(1+\frac{2(p_j-1)}{(p_j-2)(p_j^s-1)}\right)}{\prod_{p|m}\left(1+\frac{2}{p^s-1}\right)}\tilde\kappa(m,s)\Biggr\}
\\\numberthis\label{eq:deux_exp_H}&=:\frac{8^s+2^s+2}{8^s+4^s}\frac{1}{\zeta(2s)}h_m(s),
\end{align*}
with $\tilde\kappa(m,s)=1$ if $m$ is odd and equals $\frac{(-1)^{e_1}}2\frac{8^s+4^s}{8^s+2^s+2}\left(1+\frac{2(1+(-1)^{e_1})}{4^s(2^s-1)}\right)$ otherwise, and $h_m(s)$ is implicitly defined by the last equality above. By definition, $h_m(s)$ is multiplicative for any fixed real $s>1/2$. It also happens to be the expectation of $\XX_s(m)$ (which is the reason for its definition).

\begin{rema}
These quantities may seem to be artificial at first sight. The case $s=1$ is the case corresponding to our random model. The case $s=\infty$ is actually a better random model when $\alpha$ is very small compared to $x$ ($\alpha\leq 1/(2\log x)$, say), as future computations will show. Finally, the other cases are indeed artificial in the sense that the aim for their introduction is not to model any behaviour, but rather to use the expectation formalism to greatly simplify several computations to come. In a few words it will allow us to express a multiplicative function as the expectation of a totally multiplicative random variable, which will come in handy for factorizing some series as Euler products.
\end{rema}

\begin{lemm}\label{lemm:esp}
Let $m=2^{e_1} p_2^{e_2}...p_r^{e_r}$ be the prime factorization of $m$ ($e_1\geq 0$, $e_j\geq 1$ if $j=2,...,r$). Let $m_0$ be the squarefree part of $m/2^{e_1}$. Then for any real $s>1/2$:
\[\EE(\XX_s(m))=h_m(s).\]
Furthermore,
\[\EE(\XX_\infty(m))=\frac{(-1)^{\omega(m_0)}}{m_0}\prod_{\substack{2\leq j\leq r\\ e_j\text{ even}}}\left(1-\frac 2{p_j}\right)\frac{a(m)}4.\]
\end{lemm}
\begin{proo}
By the independence of the $\XX_s(p)$'s, we know that $\EE(\XX_s(\cdot))$ is a multiplicative function. Therefore, $h_\cdot(s)$ being multiplicative as well, we may check the equality only for integers $p^k$, $p$ prime and $k\geq 1$. If $p$ is odd, then if $k$ is odd, 
\[\EE(\XX_s(p^k))=a_p(s)-b_p(s)=\frac{1}p\left(\frac{2}{p^s+1}-1\right)=h_{p^k}(s).\]
If $k$ is even, 
\[\EE(\XX_s(p^k))=a_p(s)+b_p(s)=\frac{1}p\left(p-2+\frac{2}{p^s+1}\right)=h_{p^k}(s).\]
Now if $p=2$ and $k$ is odd, then
\[\EE(\XX_s(2^k))=a_2(s)-b_2(s)=\frac12\frac{4^{s}-8^{s}}{8^s+2^s+2}=h_{2^k}(s).\]
Finally, if $k$ is even, 
\[\EE(\XX_s(2^k))=a_2(s)+b_2(s)=\frac12\frac{8^{s}-4^{s}+4}{8^s+2^s+2}=h_{2^k}(s).\]
By continuity, we find the result for $\XX_\infty$ by letting $s\to\infty$, since $\tilde \kappa(m,s)\underset{s\to\infty}{\to} a(m)/4$.
\end{proo}
\begin{lemm}\label{majoration_MT}
Let $m\geq1$ and let $m_0$ be the squarefree part of $m$. We have
\[h_m(1), h_m'(1), h_m''(1)\ll \log^2(m+1)/m_0.\]
If we fix $1/2<\sigma<1$, then we also have the more general estimate
\[\sup_{\Re(s)>\sigma}|h_m(s)|\ll m^{\O(1/\log(1+m)^{\sigma/2})}/m_0.\]
\end{lemm}
\begin{proo}
First, it is clear that $\tilde\kappa(m,1)\ll 1$, and this holds true for the first and second derivatives of $\tilde\kappa(m,s)$ at $s=1$. We start by studying the product $\prod_{\substack{2\leq j\leq r\\ e_j~even}}\left(1-\frac2{p_j}\right)\left(1+\frac{2(p_j-1)}{(p_j-2)(p_j^s-1)}\right)$, and we denote the factors of this product by $u_{p_j}(s)$. We have
\[u_{p_j}(1)=1,\qquad u'_{p_j}(1)\ll \log(p_j)/p_j,\qquad u''_{p_j}(1)\ll(\log p_j)^2/p_j.\]
Thus
\[\left(\prod_{\substack{2\leq j\leq r\\ e_j~even}}u_{p_j}(s)\right)'\biggr|_{s=1}\ll \sum_{p|m} |u'_{p}(1)|\ll \log(1+m).\]
Similarly, 
\[\left(\prod_{\substack{2\leq j\leq r\\ e_j~even}}u_{p_j}(s)\right)''\biggr|_{s=1}\ll \left(\sum_{p|m}|u'_p(1)|\right)^2 +\sum_{p|m}|u''_p(1)|\ll\log^2(1+m).\]
Now we turn to the product $\prod_{p|m}\left(1+\frac{2}{p^s-1}\right)^{-1}$, which is $\ll 1$ if $s=1$, and as done above, its first derivative at $s=1$ is $\ll \log(1+m)$ and its second derivative at $s=1$ is $\ll \log^2(m+1)$. This shows that $h_m(s)\ll\log(m+1)^2/m_0$ if $s=1$, and this holds true for its first and second derivatives at $s=1$, which concludes the first part of the proof.

Let $1/2<\sigma<1$. Observe that for any $s$ with real part $>\sigma$, we have 
\[|h_m(s)|\ll \frac1{m_0}\prod_{p|m}\left(1+\O\left(\frac{1}{p^\sigma}\right)\right)\ll \frac{1}{m_0}\exp\left(\O\left(\sum_{p|m}\frac{1}{p^\sigma}\right)\right).\]
The inner sum is, if $P$ is the $\omega(m)$-th prime number:
\[\ll\sum_{p\leq P}\frac{1}{p^\sigma}\ll P^{1-\sigma}\ll \log(1+m)^{1-\sigma/2},\]
which comes from $P\ll \omega(m)\log(\omega(m))$ by the prime number theorem, and using $\omega(m)\ll \log(m+1)$. This proves that 
\[|h_m(s)|\ll m^{\O(1/\log(1+m)^{\sigma/2})}/m_0,\]
which completes the proof.
\end{proo}

Now we recall standard bounds (also found in \cite{dahl_distribution_2016}) for $d_z$, the multiplicative function defined on prime powers by $d_z(p^\nu)=\Gamma(z+\nu)/(\Gamma(z)\nu!)$. Observe that $|d_z(n)|\leq d_{|z|}(n)\leq d_k(n)$, for any $k\geq |z|$, and $d_k(mn)\leq d_k(m)d_k(n)$ for any positive integers $k,m,n$. Also note that for $k\in\NN$ and $y>3$, we have that 
\[d_k(n)e^{-n/y}\leq e^{k/y}\sum_{a_1...a_k=n}e^{-(a_1+...+a_k)/y}.\]
Therefore, 
\[\numberthis\label{eq:dz}\sum_{n=1}^\infty \frac{d_k(n)}ne^{-n/y}\leq \left(e^{1/y}\sum_{a=1}^\infty\frac{e^{-a/y}}a\right)^k\leq (\log( 3y))^k.\]

\begin{lemm}\label{lemm:majoration_sans_m}
Let $y>3$. Then for $1\leq k \leq \log(y)/3$:
\[\sum_{m=1}^\infty d_k(m)e^{-m/y}\ll y(\log (3y))^{k+2}.\]
\end{lemm}
\begin{proo}
We split the sum $\sum_{m=1}^\infty d_k(m)e^{-m/y}$ into two sums: the sum $\Sigma_1$ over $m\leq y\log^2y$, and the sum $\Sigma_2$ over $m>y\log^2 y$. By \eqref{eq:dz}, we find that 
\[\Sigma_1\leq\sum_{m\leq y\log^2y}\frac{ y\log^2y}{m}d_k(m)e^{-m/y}\ll y(\log (3y))^{k+2}.\]
Again as in \eqref{eq:dz}, the second sum is bounded by 
\begin{align*}
\Sigma_2\leq \exp\left(-\frac{\log^2 y}{2}\right)\sum_{m=1}^{\infty}d_k(m)e^{-m/(2y)}&\leq \exp\left(-\frac{\log^2 y}{2}\right) \left(e^{1/(2y)}\sum_{m=1}^\infty e^{-m/(2y)}\right)^k
\\&\ll\exp\left(-\frac{\log^2 y}{2}\right)(2y)^k\ll 1.
\end{align*}
This proves the result.
\end{proo}
\begin{lemm}\label{lemm:majoration-m0}
Let $x$ be large, $(a_m)$ be a complex sequence such that $a_m\ll m^{\O(1/\log(1+m)^{\kappa}+1/\log_2 x)}/m_0$, where $m_0$ is the squarefree part of $m$ and $0<\kappa<1/2$ is a real number. Then for $B>0$, $z\in\CC$ such that $|z|\leq B\log x/(\log_2 x\log_3 x)$:
\[\sum_{m=1}^\infty\frac{d_z(m)}ma_m\ll \exp\left((2B+o(1))\frac{\log x}{\log_2 x}\right).\]
\end{lemm}
\begin{proo}
We split the above sum into two sums: one over $m\leq e^{(\log\log x)^{1/\kappa}}$ and the other one over $m>e^{(\log\log x)^{1/\kappa}}$. We let $m=m_0m_1^2$ and put $k=\lceil |z|\rceil$. Since $d_k(m)\leq d_k(m_0)d_k(m_1^2)\leq d_k(m_0)d_k(m_1)^2$, the first sum is 
\begin{align*}
&\ll \sum_{m=1}^\infty \frac{d_k(m)\exp\left(\O\left(\log_2^{1/\kappa-1}x\right)\right)}{m_0m}\ll \exp\left(\O\left(\log_2^{1/\kappa-1}x\right)\right)\sum_{m= 1}^\infty\frac{d_k(m)}{m_0m}
\\&\ll \exp\left(\O\left(\log_2^{1/\kappa-1}x\right)\right)\sum_{m_0=1}^\infty\frac{d_k(m_0)}{m_0^2}\sum_{m_1=1}^\infty\frac{d_k(m_1)^2}{m_1^2}
=\exp\left(\O\left(\log_2^{1/\kappa-1}x\right)\right)\zeta(2)^k\sum_{m=1}^\infty\frac{d_k(m)^2}{m^2}.
\end{align*}
Lemma 3.3 of \cite{lamzouri_extreme_2010} reveals that this expression is 
\[\ll\exp\left(\O\left(\log_2^{1/\kappa-1}x\right)\right)\zeta(2)^k\exp((2+o(1))k\log_2(k+3))\ll \exp\left((2B+o(1))\frac{\log x}{\log_2 x}\right),\]
Similarly the other sum is
\begin{align*}
&\ll\sum_{m>e^{(\log\log x)^{1/\kappa}}}\frac{d_z(m)m^{\O(1/\log(1+m)^{\kappa}+1/\log_2 x)}}{mm_0}
\ll\sum_{m=1}^\infty\frac{d_k(m)m^{\O(1/\log\log x)}}{mm_0}
\\&\ll \sum_{m_0=1}^\infty \frac{d_k(m_0)}{m_0^{2-\O(1/\log\log x)}}\sum_{m_1=1}^\infty\frac{d_k(m_1)^2}{m_1^{2-\O(1/\log\log x)}}.
\end{align*}
Once again, Lemma 3.3 of \cite{lamzouri_extreme_2010} shows that this expression above is $\ll\exp\left((2B+o(1))\frac{\log x}{\log_2 x}\right)$, which concludes the proof.
\end{proo}

\begin{lemm}\label{lemm:prolongement} 
Let $x$ be large, $B>0$ be a constant and $z\in \CC$ be such that $|z|\leq B\log x/(\log_2 x\log_3 x)$, let $s$ be a complex with real part $>1/2$. Then
\[\sum_{m=1}^\infty\frac{d_z(m)}mh_m(s)=\prod_p\left(a_p(s)\left(1-\frac1p\right)^{-z}+b_p(s)\left(1+\frac1p\right)^{-z}+c_p(s)\right)\]
and both of these terms are analytic in the variable $s$ in this region.
\end{lemm}
\begin{proo}
This equality is true for any real number $s>1/2$, since in this case \nameref{lemm:esp} ensures that
\[\sum_{m=1}^\infty\frac{d_z(m)}mh_m(s)=\sum_{m=1}^\infty\frac{d_z(m)}m\EE(\XX_s(m))=\EE(L(1,\XX_s)^z),\]
which concludes (in the case of real numbers), thanks to the Euler product expression for $L(1,\XX_s)$. To prove that the equality holds for complex numbers, we shall rely on the analytic continuation principle. For any $1/2<\sigma<1$ and any $s$ such that $\Re(s)>\sigma$, \nameref{majoration_MT} and \nameref{lemm:majoration-m0} ensure that
\[\sum_{m=1}^\infty\sup_{s:\Re(s)>\sigma}\left|\frac{d_z(m)}mh_m(s)\right|\ll \sum_{m=1}^\infty \frac{d_{|z|}(m)}{m_0m}m^{\O(1/\log(1+m)^{\sigma/2})}\ll_{\sigma,x}1.\]
This proves that $\sum_{m=1}^\infty\frac{d_z(m)}mh_m(s)$ is a series of analytic functions that converges uniformly in any compact of the half-plane $\Re(s)>1/2$, and hence is analytic in this region. Furthermore in this region, we have $a_p(s),b_p(s)=1/2+\O(1/p)$, which leads to\footnote{Recall that $c_p(s)=1-a_p(s)-b_p(s)$.}
\begin{align*}
&\prod_{p\geq 3}\left(a_p(s)\left(1-\frac1p\right)^{-z}+b_p(s)\left(1+\frac1p\right)^{-z}+c_p(s)\right)
\\&=\prod_{p\geq 3}\left(a_p(s)\left[\left(1-\frac1p\right)^{-z}-1\right]+b_p(s)\left[\left(1+\frac1p\right)^{-z}-1\right]+1\right)
\\&=\prod_{p\geq 3}\left(\left[\frac12+\O\left(\frac1p\right)\right]\left[\frac zp+\O_z\left(\frac1{p^2}\right)\right]+\left[\frac12+\O\left(\frac1p\right)\right]\left[-\frac zp+\O_z\left(\frac1{p^2}\right)\right]+1\right)
\\&=\prod_{p\geq 3}\left(1+\O_z\left(\frac 1{p^2}\right)\right),
\end{align*}
which is analytic in $\Re(s)>1/2$. The factor corresponding to $p=2$ is also analytic in this region. This concludes the proof.
\end{proo}
\begin{rema}
We may sometimes write $\EE(L(1,\XX_s)^z)$ even when $s$ is a complex number. This should not be defined, since $(\XX_s(n))$ are not defined if $s\not\in\RR$. However, thanks to the result above, we can extend the definitions and denote the Euler product of \nameref{lemm:prolongement} above by $\EE(L(1,\XX_s)^z)$.
\end{rema}
\begin{lemm}\label{lemm:exp_phi}
Let $z\in\CC$. Then
\[\sum_{m=1}^\infty\frac{d_z(m)}mH_m(1)=\frac{\EE(L(1,\XX)^z)}{\zeta(2)}\]
and
\[\sum_{m=1}^\infty\frac{d_z(m)}mH'_m(1)=\frac{\EE(L(1,\XX)^z)}{\zeta(2)}(\phi(z)-2\gamma),\]
with $\phi$ being defined in \nameref{theo-hd}.
\end{lemm}
\begin{proo}
Let $s>1/2$ be a real number and put $\psi(s):=\frac{8^s+2^s+2}{8^s+4^s}$ for commodity. By \eqref{eq:deux_exp_H} and by \nameref{lemm:esp}:
\begin{align*}
\sum_{m=1}^\infty\frac{d_z(m)}mH_m(s)&=\frac{\psi(s)}{\zeta(2s)}\sum_{m=1}^\infty\frac{d_z(m)}mh_m(s)=\frac{\psi(s)}{\zeta(2s)}\EE\left(\sum_{m=1}^\infty\frac{d_z(m)}m\XX_s(m)\right)=\frac{\psi(s)}{\zeta(2s)}\EE(L(1,\XX_s)^z).
\end{align*}
At $s=1$, we simply find the expected $\EE(L(1,\XX)^z)/\zeta(2)$. Furthermore,
\begin{align*}
\sum_{m=1}^\infty\frac{d_z(m)}mH'_m(1)&=\frac{\psi'(1)}{\zeta(2)}\EE(L(1,\XX)^z)-\frac{2\zeta'(2)\psi(1)}{\zeta(2)^2}\EE(L(1,\XX)^z)+\frac{\psi(1)}{\zeta(2)}\frac{\diff}{\diff s}\EE(L(1,\XX_s)^z)\bigr|_{s=1}
\\\numberthis\label{eq:lemme_phi}&=\frac{\EE(L(1,\XX)^z)}{\zeta(2)}\left\{-\frac{\log 2}2-\frac{2\zeta'(2)}{\zeta(2)}+\frac{\frac{\diff}{\diff s}\EE(L(1,\XX_s)^z)\bigr|_{s=1}}{\EE(L(1,\XX)^z)}\right\}.
\end{align*}

Recalling the definitions \eqref{eq:def_probas} of $a_p(s)$, $b_p(s)$ and that $c_p(s)=1-a_p(s)-b_p(s)$, we may use the Euler product expression for $\EE(L(1,\XX_s)^z)$ to find that
\begin{align*}
\frac{\frac{\diff}{\diff s}\EE(L(1,\XX_s)^z)\bigr|_{s=1}}{\EE(L(1,\XX)^z)}&=\frac{\diff}{\diff s}\log \EE(L(1,\XX_s)^z)\biggr|_{s=1}
\\&=\sum_{p\geq 3} \frac{a'_p(1)\left[\left(1-\frac1p\right)^{-z}-1\right]}{a_p(1)\left(1-\frac1p\right)^{-z}+b_p(1)\left(1+\frac1p\right)^{-z}+c_p(1)}+\frac{\log 2}9\frac{\frac{-13}42^z+\frac94\left(\frac{2}3\right)^z+1}{\frac16 2^z+\frac 12\left(\frac23\right)^z+\frac43}
\\\numberthis\label{eq:exp_phi_3}&=\sum_{p\geq 3} \frac{\frac{2\log p}{(p+1)^2}\left[1-\left(1-\frac1p\right)^{-z}\right]}{\EE\left(\left(1-\frac{\XX(p)}p\right)^z\right)}+\frac{\log 2}9\frac{\frac{-13}42^z+\frac94\left(\frac{2}3\right)^z+1}{\frac16 2^z+\frac 12\left(\frac23\right)^z+\frac43}.
\end{align*}
In \eqref{eq:lemme_phi}, this proves that 
\[\sum_{m=1}^\infty\frac{d_z(m)}mH'_m(1)=\frac{\EE(L(1,\XX)^z)}{\zeta(2)}(\phi(z)-2\gamma),\]
which is the expected result.
\end{proo}
\begin{lemm}\label{majoration_esp_s}
Let $z\in\CC$, $\Re(z)\geq -1$. Then for any $s$ with real part $\geq 2/3$:
\[\EE(L(1,\XX_s)^z)\ll \EE(L(1,\XX)^{\Re(z)}).\]
This also holds for $\XX_\infty$.

Furthermore, this estimate above holds if the left-hand side is replaced by $\frac{\diff}{\diff s} \EE(L(1,\XX_s)^z)\biggr|_{s=1}$ or $\frac{\diff^2}{\diff s^2} \EE(L(1,\XX_s)^z)\biggr|_{s=1}$.
\end{lemm}
\begin{proo}
We put $\sigma:=\Re(z)$. First, we assume that $-1\leq \sigma\leq 8$. In this case for $p\geq 3$,
\begin{align*}
&\left|a_p(s)\left(1-\frac1p\right)^{-z}+b_p(s)\left(1+\frac1p\right)^{-z}+c_p(s)\right|
\\&\leq \left[\frac12-\frac3{2p}+\O(p^{-5/3})\right]\left(1-\frac 1p\right)^{-\sigma}+\left[\frac12-\frac{1}{2p}\right]\left(1+\frac1p\right)^{-\sigma}+\frac2p+\O(p^{-5/3})
\\&=\left[\frac12-\frac3{2p}+\O(p^{-5/3})\right]\left\{\left(1-\frac 1p\right)^{-\sigma}-1\right\}+\left[\frac12-\frac{1}{2p}\right]\left\{\left(1+\frac1p\right)^{-\sigma}-1\right\}+1+\O(p^{-5/3})
\\&=1+\O(p^{-5/3}),
\end{align*}
which is obtained by Taylor expanding the factors $(1\pm 1/p)^{-\sigma}$. Taking the product over primes $\geq 3$ and using the fact that the factor corresponding to $p=2$ is $\O(1)$, we find that 
\[\EE(L(1,\XX_s)^z)\ll \prod_p(1+\O(p^{-5/3}))\ll 1 \ll \EE(L(1,\XX)^{\sigma}),\]
which is the expected result.
 
If $\sigma>8$, notice that for any $p\geq 2$:
\begin{align*}
\left|a_p (1)\left(1-\frac1p\right)^{-z}+b_p(1)\left(1+\frac1p\right)^{-z}+c_p(1)\right|&\geq a_p(1) \left(1-\frac1p\right)^{-8}-b_p(1)\left(1+\frac1p\right)^{-8}-c_p(1)
\\&\numberthis\label{eq:mino_1}\geq 1>0,
\end{align*}
and hence the following computations are licit\footnote{Observe that $b_p(s)=b_p(1)$.}:
\begin{align*}
\frac{\EE(L(1,\XX_s)^z)}{\EE(L(1,\XX)^z)}&\ll \frac{\prod_{p\geq 3}\left|a_p(s)\left[\left(1-\frac1p\right)^{-z}-1\right]+b_p(1)\left[\left(1+\frac1p\right)^{-z}-1\right]+1\right|}{\prod_{p\geq 3}\left|a_p(1)\left[\left(1-\frac1p\right)^{-z}-1\right]+b_p(1)\left[\left(1+\frac1p\right)^{-z}-1\right]+1\right|}
\\\numberthis\label{eq:inter_prod}&=\prod_{p\geq 3}\left|1+\frac{(a_p(s)-a_p(1))\left[\left(1-\frac1p\right)^{-z}-1\right]}{a_p(1)\left[\left(1-\frac1p\right)^{-z}-1\right]+b_p(1)\left[\left(1+\frac1p\right)^{-z}-1\right]+1}\right|.
\end{align*}
By definition, since $\Re(s)\geq 2/3$:
\[\numberthis\label{eq:diff_a} a_p(s)-a_p(1)\ll 1/p^{5/3}.\]
If $\left(1-\frac1p\right)^{-\sigma}> e^{10}$, then we see that
\begin{align*}
\left|a_p (1)\left[\left(1-\frac1p\right)^{-z}-1\right]+b_p(1)\left[\left(1+\frac1p\right)^{-z}-1\right]+1\right|&\geq  \frac1{10}\left[\left(1-\frac1p\right)^{-\sigma}-1\right],
\end{align*}
and if $\left(1-\frac1p\right)^{-\sigma}\leq e^{10}$, then we may use \eqref{eq:mino_1} to show that the denominator in \eqref{eq:inter_prod} is $\geq 1$ and \eqref{eq:diff_a} to show that the numerator is $\ll p^{-5/3}$. Either way,
\[\frac{\EE(L(1,\XX_s)^z)}{\EE(L(1,\XX)^z)}\ll\prod_{p\geq 3}\left(1+\O\left(\frac{1}{p^{5/3}}\right)\right)\ll 1.\]
Since $\EE(L(1,\XX)^z)\ll\EE(L(1,\XX)^\sigma)$, this implies that for any value of $\sigma\geq -1$, we have 
\[\EE(L(1,\XX_s)^z)\ll\EE(L(1,\XX)^\sigma).\]
By the analyticity proved in \nameref{lemm:prolongement} above, we may let $s\to \infty$, which concludes the first part of the proof.

We now turn to the other two estimates. We first bound $\frac{\diff}{\diff s} \EE(L(1,\XX_s)^z)\biggr|_{s=1}$. By our \eqref{eq:exp_phi_3} above, we know that 
\begin{align*}\numberthis\label{eq:exp_inter}
\frac{\diff}{\diff s}\EE(L(1,\XX_s)^z)\bigr|_{s=1}&=\EE(L(1,\XX)^z)\left(\sum_{p\geq 3} \frac{\frac{2\log p}{(p+1)^2}\left[1-\left(1-\frac1p\right)^{-z}\right]}{\EE\left(\left(1-\frac{\XX(p)}p\right)^z\right)}+\frac{\log 2}9\frac{\frac{-13}42^z+\frac94\left(\frac{2}3\right)^z+1}{\frac16 2^z+\frac 12\left(\frac23\right)^z+\frac43}\right).
\end{align*}
Observe that, for a given prime $p$ and $-1\leq \sigma\leq 8$, we have $\prod_{q\neq p}\EE\left(\left(1-\frac{\XX(q)}q\right)^{-z}\right)\ll \EE(L(1,\XX)^\sigma)$, since the missing factor $\EE((1-\XX(p)/p)^{-\sigma})$ is $\gg 1$. Thus following the same reasoning as above, if $-1\leq \sigma\leq 8$ we have 
\begin{align*}
\frac{\diff}{\diff s}\EE(L(1,\XX_s)^z)\bigr|_{s=1}&\ll \EE(L(1,\XX)^\sigma)\left(\sum_{p\geq 3} \frac{2\log p}{(p+1)^2}\left[1+\left(1-\frac1p\right)^{-\sigma}\right]+\O(1)\right)
\\&\ll \EE(L(1,\XX)^\sigma).
\end{align*}
Otherwise if $\sigma>8$, then by considering the cases $\left(1-\frac1p\right)^{-\sigma}> e^{10}$ and $\left(1-\frac1p\right)^{-\sigma}\leq e^{10}$, we deduce that the term inside parenthesis in \eqref{eq:exp_inter} is $\ll 1$. Therefore, for all $\sigma\geq -1$, we have 
\[\frac{\diff}{\diff s}\EE(L(1,\XX_s)^z)\bigr|_{s=1}\ll\EE(L(1,\XX)^\sigma).\]

We now turn to the quantity $\frac{\diff^2}{\diff s^2} \EE(L(1,\XX_s)^z)\biggr|_{s=1}$. For convenience, we put $E_{z,p}(s):=\EE\left(\left(1-\frac{\XX_s(p)}p\right)^{-z}\right)$, so that $\EE(L(1,\XX_s)^z)=\prod_p E_{z,p}(s)$. We have
\begin{align*}
\numberthis\label{eq:X1}\frac{\diff^2}{\diff s^2} \EE(L(1,\XX_s)^z)\biggr|_{s=1}=\sum_{p}E_{z,p}'(1)\sum_{q\neq p}\left[\prod_{r\neq p,q}E_{z,r}(1)\right] E_{z,q}'(1)+\sum_p E_{z,p}''(1)\prod_{q\neq p} E_{z,q}(1).
\end{align*}
Now and for the rest of the proof, observe that we have for $p\geq 3$:
\[\numberthis\label{eq:X2}E''_{z,p}(1)= \frac{2(p-1)(\log p)^2}{(p+1)^3}\left(\left(1-\frac1p\right)^{-z}-1\right),\]
and recall that in \eqref{eq:exp_inter} we showed that 
\[\numberthis\label{eq:X3} E'_{z,p}(1)=\frac{2\log p}{(p+1)^2}\left[1-\left(1-\frac1p\right)^{-z}\right] .\]
Once again, we do the exact same reasoning: if $-1\leq \sigma\leq 8$, then we can introduce the missing factors in the products of \eqref{eq:X1} (respectively $E_{\sigma,p}(1)\times E_{\sigma,q}(1)$ and $E_{\sigma,p}(1)$, which are $\gg 1$), and hence 
\[\prod_{r\neq p,q}E_{z,r}(1),\prod_{q\neq p} E_{z,q}(1)\ll \prod_p E_{z,p}(1)=\EE(L(1,\XX)^z).\]
Therefore
\[\frac{\diff^2}{\diff s^2} \EE(L(1,\XX_s)^z)\biggr|_{s=1}\ll \EE(L(1,\XX)^\sigma)\left[\left(\sum_p |E_{z,p}'(1)|\right)^2+\sum_p |E_{z,p}''(1)|\right]\ll \EE(L(1,\XX)^\sigma).\]

Now if $\sigma> 8$ then $|E_{z,p}(1)|\geq 1$ by \eqref{eq:mino_1} and hence it is licit to deduce from \eqref{eq:X1} that 
\[\frac{\diff^2}{\diff s^2} \EE(L(1,\XX_s)^z)\biggr|_{s=1}\ll \EE(L(1,\XX)^\sigma)\left[\left(\sum_p \left|\frac{E_{z,p}'(1)}{E_{z,p}(1)}\right|\right)^2+\sum_p \left|\frac{E_{z,p}''(1)}{E_{z,p}(1)}\right|\right].\]
As shown earlier when bounding the term inside the parenthesis in \eqref{eq:exp_inter}, by considering the cases $\left(1-\frac1p\right)^{-\sigma}> e^{10}$ and $\left(1-\frac1p\right)^{-\sigma}\leq e^{10}$ and using \eqref{eq:X3} we show that the first sum above is $\ll 1$ (the case $p=2$ is dealt with separately). Now we turn to the second sum and once again considering the cases $\left(1-\frac1p\right)^{-\sigma}> e^{10}$ and $\left(1-\frac1p\right)^{-\sigma}\leq e^{10}$ and using \eqref{eq:X2}, this second sum is $\ll 1$, and hence for all $\sigma\geq -1$,
\[\frac{\diff^2}{\diff s^2} \EE(L(1,\XX_s)^z)\biggr|_{s=1}\ll \EE(L(1,\XX)^{\sigma}).\]
\end{proo}

Observe that, as in \eqref{eq:lemme_phi}, if $f(m)$ is a linear combination of $H_m(1)$, $H'_m(1)$ and $H''_m(1)$, then 
\[\sum_{m=1}^\infty \frac{d_z(m)}mf(m)=A\EE(L(1,\XX)^z)+B\frac{\diff}{\diff s}\EE(L(1,\XX_s)^z)\biggr|_{s=1}+C\frac{\diff^2}{\diff s^2}\EE(L(1,\XX_s)^z)\biggr|_{s=1}\]
for some constants $A$, $B$, and $C$. Therefore, our \nameref{majoration_esp_s} implies the following corollary:
\begin{coro}\label{coro:majo_esp}
Let $f(m)$ be a linear combination of $H_m(1)$, $H'_m(1)$ and $H''_m(1)$. Let $z$ be a complex number with real part $\geq -1$. Then
\[\sum_{m=1}^\infty \frac{d_z(m)}mf(m)\ll \EE(L(1,\XX)^{\Re(z)}).\]
\end{coro}
\begin{rema}\label{rema:phi}
In particular, thanks to \nameref{lemm:exp_phi} combined to the corollary above, we know that 
\[\EE(L(1,\XX)^z)\phi(z)\ll \EE(L(1,\XX)^{\Re(z)}).\]
\end{rema}

\section{The sum \texorpdfstring{$\sum_{d\in\D_\alpha(x)}\chi_d(m)$}{of characters}}
We recall that $X_\alpha:=x^\alpha-x^{-1-\alpha}$, $\alpha_0$ is such that $X_{\alpha_0}=1$ and $\alpha_1$ is such that $X_{\alpha_1}=2$, as defined above \eqref{eq:def_alpha}.
\begin{prop}\label{prop:sum-chid}
Let $m$ be a positive integer, $x$ be a large real number. We recall that $Y_1$ is defined in \nameref{lemm}. Then uniformly in $\alpha_1\leq\alpha<1/2$, we have
\begin{align*}
\sum_{d\in\D_\alpha(x)}\chi_d(m)&=\sqrt x\left(\alpha^2\frac{\EE(\XX(m))}{2\zeta(2)}\log^2x+\alpha\{H'_m(1)+(2\gamma-2\alpha) \frac{\EE(\XX(m))}{\zeta(2)}\}(\log x-2)+\eth_2(m)\right)
\\&+(1-\sqrt{Y_1(1,\alpha)})\EE(\XX_\infty(m))+\mathcal I(x,m,\alpha)+ \O\left(m x^{\alpha}\log x+\frac{\log^2(m+1)}{m_0}\right)
\end{align*}
where 
\[\eth_2:=\eth_2(m):=\frac{H''_m(1)}2+2\gamma H_m'(1)+(\gamma^2-2\gamma_1)H_m(1)\] 
with $\gamma_1$ being the first Stieltjes constant, and
\begin{align*}
\mathcal I(x,m,\alpha)&:=\frac{X_\alpha^{1/(2\alpha)}}{2\pi i}\int_{-1/\log\log x-i\infty}^{-1/\log\log x+i \infty}H_m(1+s)\zeta^2(1+s)\frac{1/(2\alpha)}{s(s+1/(2\alpha))}X_\alpha^s \diff s
\\&\ll \sqrt x \frac{m^{1/\log(1+m)^{1/4}}}{m_0}(\log\log x)^3.
\end{align*}
Uniformly in $\alpha_0\leq \alpha\leq\alpha_1$, we have 
\[\sum_{d\in\D_\alpha(x)}\chi_d(m)=(\sqrt x-\sqrt{Y_1(1,\alpha)})\EE(\XX_\infty(m))+\O\left(m \log x\right).\]
\end{prop}
\begin{rema}
This proves that for any fixed $0<\alpha<1/2$:
\[\sum_{d\in\D_\alpha(x)}\chi_d(m)=\alpha^2\sqrt x\frac{\EE(\XX(m))}{2\zeta(2)}\log^2x+\O_m(\sqrt x\log x),\]
and hence 
\[\numberthis\label{eq:card_D}\#\D_\alpha(x)\sim \frac{\alpha^2\sqrt x}{2\zeta(2)}\log^2x.\]
Therefore, for any fixed $m$ and fixed $0<\alpha<1/2$, 
\[\frac{1}{\EE(\XX(m))}\sum_{d\in\D_\alpha(x)}\chi_d(m)\sim \#\D_\alpha(x),\]
which is the very reason for the definition of $\XX(m)$'s law. However we chose not to present our proposition in this way, for the error terms we would get would not be sharp enough for what the next section requires.
\end{rema}
\begin{proo}[of \nameref{prop:sum-chid}]
We recall the equality proved in \eqref{main_term} and the definition of $B_m(u)$ in \nameref{lemm:def_bm}, so that using \nameref{esti_ca} leads to
\begin{align*}
\sum_{d\in\D_{\alpha}(x)}\chi_d(m)&=\sum_{u=1}^{X_\alpha}\frac{\sqrt x-\sqrt{Y_1}}{4u}\sum_{\substack{2<a\leq 4u^2+2\\d(a,u)\in\D}}\frac{C_{m,a,u}}m+\O\left(m\sum_{u=1}^{X_\alpha}\rho(u)\right)
\\&=\sum_{u=1}^{X_\alpha}\frac{\sqrt x-\sqrt{Y_1}}{4u}\frac{(-1)^{\omega(m_0)}}{m_0}\prod_{\substack{2\leq j\leq r\\ e_j \text{ is even}}}\left(1-\frac 2{p_j}\right)\prod_{\substack{p_j|u\\ e_j\text{ is even}}}\left(1+\frac{1}{p_j-2}\right)B_m(u)
\\&\numberthis\label{mainterm2}\qquad +\O(mx^\alpha\log x),
\end{align*}
since $\rho(u)\ll d(u)$, with $d(u)$ being the usual divisor function, by the Chinese remainder Theorem.

Now if 
\[\mathfrak F_m(u):=\frac{(-1)^{\omega(m_0)}}{4m_0}\prod_{\substack{2\leq j\leq r\\ e_j\text{ is even}}}\left(1-\frac2{p_j}\right)\prod_{\substack{p_j|u\\ e_j\text{ is even}}}\left(1+\frac{1}{p_j-2}\right)B_m(u),\]
\nameref{lemm:function_perron} reveals that for $\Re(s)>1$:
\begin{align*}
\L_m(s):=\sum_{\substack{u\geq 1\\(u,m_0)=1}}\frac{\mathfrak F_m(u)}{u^s}&=\frac{(-1)^{\omega(m_0)}}{m_0}\prod_{\substack{2\leq j\leq r\\e_j\text{ is even}}}\left(1-\frac2{p_j}\right)\frac{a(m)}4\times \Sigma(m,s)
\\&=\frac{(-1)^{\omega(m_0)}}{m_0}\prod_{\substack{2\leq j\leq r\\ e_j\text{ is even}}}\left(1-\frac2{p_j}\right)\frac{a(m)}4\times G_m(s)\zeta^2(s)=H_m(s)\zeta(s)^2
\end{align*}
where $H_m(s)$ was defined in \eqref{defH}. Observe that $H_m(s)$ is analytic in the region $\Re(s)>1/2$, since $G_m(s)$ is. From now on, we assume that $\alpha\geq \alpha_1$, from which we deduce that $\alpha\gg 1/\log x$ and hence 
\[\numberthis\label{eq:devAsympX} X_\alpha^{1/(2\alpha)}=\sqrt x(1-x^{-1-2\alpha})^{1/(2\alpha)}=\sqrt x\exp\left(\O\left(\frac{1}{\alpha x}\right)\right)=\sqrt x+\O\left(1/x^{1/3} \right).\]
Therefore, using \nameref{lemm:size_Y}, \eqref{eq:devAsympX} and the fact that $\frak F_m(u)\ll \log(1+u)/m_0$, we may deduce that 
\begin{align*}
&\sum_{u\leq X_\alpha}(\sqrt x-\sqrt{Y_1})\frac{\frak F_m(u)}u
\\&=\sum_{u\leq X_\alpha}(X_\alpha^{1/(2\alpha)}-u^{1/(2\alpha)})\frac{\frak F_m(u)}u+\sum_{u\leq X_\alpha}(\sqrt x-X_\alpha^{1/(2\alpha)})\frac{\frak F_m(u)}u+\sum_{u\leq X_\alpha}(u^{1/(2\alpha)}-\sqrt{Y_1(u,\alpha)})\frac{\frak F_m(u)}u
\\&\numberthis\label{eq:avantperron}=\sum_{u\leq X_\alpha}(X_\alpha^{1/(2\alpha)}-u^{1/(2\alpha)})\frac{\frak F_m(u)}u+(1-\sqrt{Y_1(1,\alpha)})\EE(\XX_\infty(m))+\O\left(\frac{1}{m_0}\right).
\end{align*}
The expectation appears since $\EE(\XX_\infty(m))=\mathfrak F_m(1)$ by \nameref{lemm:esp}, because\footnote{Recall that $a(m)=4$ if $m$ is odd, and equals $2(-1)^{e_1}$ otherwise.} $B_m(1)=a(m)$ by \nameref{lemm:def_bm}.

We now apply a weighted version of Perron's formula -which has the advantage of being absolutely convergent-, as a direct consequence of Theorem 2.1 of part II.2 of \cite{tenenbaum}. We fix $\kappa>0$, and
\[\numberthis\label{sum:fu}\sum_{u\leq X_\alpha}(X_\alpha^{1/(2\alpha)}-u^{1/(2\alpha)})\frac{\mathfrak F_m(u)}{u}=\frac{X_\alpha^{1/(2\alpha)}}{2\pi i}\int_{\kappa-i\infty}^{\kappa+i \infty}\L_m(1+s)\frac{1/(2\alpha)}{s(s+1/(2\alpha))}X_\alpha^s \diff s.\]

Moving the integration segment to the left\footnote{This can be done since the integrand is of finite order.}, say to $\Re(s)=-1/\log\log x$, we can apply the residue theorem to the only pole of the integrand: a third order pole at $s=0$. The said residue is equal to 
\begin{align*}
\Res=\frac{1}{2!(2\alpha)}\frac{\diff ^2}{\diff s^2}\left(\frac{H_m(1+s)\zeta^2(1+s)s^2}{s+1/(2\alpha)}X_\alpha^s\right)\Bigr|_{s\to 0}&=\frac{ H_m(1)}2\log^2 X_\alpha+\Delta_1\log X_\alpha+\Delta_2.
\end{align*}
where
\[\Delta_k= \frac{1}{k}\frac{\diff^k}{\diff s^k}\frac{H_m(1+s)\zeta^2(1+s)s^2}{2\alpha s+1}\Bigr|_{s\to 0},\quad k=1,2.\]
If we let 
\[\eth_k= \frac{1}{k}\frac{\diff^k}{\diff s^k}H_m(1+s)\zeta^2(1+s)s^2\Bigr|_{s\to 0},\quad k=1,2,\]
then, 
\[\numberthis\label{eq:defEth}\eth_1=H'_m(1)+2\gamma H_m(1),\quad \eth_2=\frac{H''_m(1)}2+2\gamma H_m'(1)+(\gamma^2-2\gamma_1)H_m(1),\]
and
\[\Delta_1=\eth_1-2\alpha H_m(1),\qquad \Delta_2=\eth_2-2\alpha\eth_1+(2\alpha)^2H_m(1).\]
Therefore, 
\[\Res=\frac{ H_m(1)}2\log^2 X_\alpha+\Bigl\{H'_m(1)+[2\gamma-2\alpha] H_m(1)\Bigr\}(\log X_\alpha-2\alpha)+\eth_2.\]

Consequently, using \eqref{sum:fu} in \eqref{eq:avantperron} leads to 
\begin{align*}
\sum_{u\leq X_\alpha}(\sqrt x-\sqrt{Y_1})\frac{\mathfrak F_m(u)}{u}=X_\alpha^{1/(2\alpha)}\Res+(1-\sqrt{Y_1(1,\alpha)})\EE(\XX_\infty(m))+\mathcal I(x,m,\alpha)+\O\left(\frac{1 }{m_0}\right),
\end{align*}
where
\[\numberthis\label{eq:defI}\mathcal I(x,m,\alpha):=\frac{X_\alpha^{1/(2\alpha)}}{2\pi i}\int_{-1/\log\log x-i\infty}^{-1/\log\log x+i \infty}H_m(1+s)\zeta^2(1+s)\frac{1/(2\alpha)}{s(s+1/(2\alpha))}X_\alpha^s \diff s.\]

Inserting this estimate in \eqref{mainterm2}, we get
\begin{align*}
\sum_{d\in\D_\alpha(x)}\chi_d(m)&=\sum_{u\leq X_\alpha}\frac{(\sqrt x-\sqrt{Y_1})\frak F_m(u)}u+ \O\left(m x^{\alpha}\log x\right)
\\&\numberthis\label{eq:sum_chi_d_final}=X_\alpha^{1/(2\alpha)}\Res+(1-\sqrt{Y_1(1,\alpha)})\EE(\XX_\infty(m))+\mathcal I(x,m,\alpha)+\O\left(m x^{\alpha}\log x+\frac{1 }{m_0}\right).
\end{align*}

Naturally we wish to use that $X_\alpha\approx x^\alpha$ to simplify the above expression. Since $X_\alpha=x^\alpha-x^{-1-\alpha}$, we may write that 
\[\numberthis\label{eq:logXa}\log X_\alpha= \log (x^\alpha)+\log(1-x^{-1-2\alpha})=\alpha \log x+\O(1/x),\]
and
\[\log^2 X_\alpha =\alpha^2\log^2 x+\O(\log x/x).\]
Also recall the estimate \eqref{eq:devAsympX}. Therefore changing these expressions in $\Res$ by the main term of the above estimates introduces an error 
\[\ll \frac{H_m(1)+H_m'(1)+H_m''(1)}{x^{1/3}}\ll \frac{h_m(1)+h_m'(1)+h_m''(1)}{x^{1/3}}\ll \frac{\log^2(m+1)}{m_0 x^{1/3}},\]
by \nameref{majoration_MT}. Thus \eqref{eq:sum_chi_d_final} becomes
\begin{align*}
\sum_{d\in\D_\alpha(x)}\chi_d(m)&=\sqrt x\left(\alpha^2\frac{H_m(1)}2\log^2x+\alpha\{H'_m(1)+(2\gamma-2\alpha) H_m(1)\}(\log x-2)+\eth_2\right)
\\&+(1-\sqrt{Y_1(1,\alpha)})\EE(\XX_\infty(m))+\mathcal I(x,m,\alpha)+ \O\left(m x^{\alpha}\log x+\frac{\log^2(m+1)}{m_0}\right).
\end{align*}
This concludes the first part of the proof, since $H_m(1)=\EE(\XX(m))/\zeta(2)$. 

To prove the bound on $\mathcal I(x,m,\alpha)$, we will rely on estimates for $\zeta(s)$ in the critical strip and more precisely the fact that for each positive constant $c$ and $t\geq 2$ such that $\sigma\geq 1-c/\log t$:
\[\numberthis\label{eq:zeta_line}\zeta(\sigma+it)\ll \log t\]
(see Theorem II.3.9 of \cite{tenenbaum} for a proof). Furthermore let $\mu(\sigma)$ be the infinimum of all real numbers $a$ such that $\zeta(\sigma+i\tau)=\O(\tau^a)$, $|\tau|\geq1$, $0<\sigma\leq 1$. It is well known that $\mu$ is convex, thanks to the Phragmén-Lindelöf theorem. It is also known that $\mu(1/2)\leq 1/6$. Many improvements of this bound have been given since, but would not bring significantly better results here. These facts may be found through the section II.3.4 of \cite{tenenbaum}. Since $\mu(\sigma)\leq \frac13(1-\sigma)$ for any $\sigma\in[1/2,1]$, we find that
\[\numberthis\label{eq:mu}\mu\left(1-\frac{1}{\log\log x}\right)\leq \varepsilon,\]
for any small and fixed $\varepsilon>0$. Consequently, since \nameref{majoration_MT} reveals that $H_m(1+s)\ll h_m(1+s)\ll m^{1/\log(1+m)^{1/4}}/m_0$ on the line $\Re(s)=-1/\log\log x$ and since $\zeta(1+s)\ll1/|s|$ when $s\to 0$, we find that 
\begin{align*}
\mathcal I(x,m,\alpha)&\ll \sqrt x \frac{m^{1/\log(1+m)^{1/4}}}{m_0}\left(\log\log^3 x+\int_2^\infty \frac{|\zeta(1-1/\log\log x+it)|^2}{t\sqrt{(2\alpha t)^2+1}}\diff t\right).
\end{align*}
The contribution to integral above when $t\geq (\log x)^2$ is relatively small. Indeed in this range, we have $\sqrt t\alpha \geq \sqrt t\alpha_1\gg \sqrt t/\log x\geq 1$, and hence $(t\alpha)^2\gg t$. With \eqref{eq:mu}, this implies that
\begin{align*}
\int_{(\log x)^2}^\infty \frac{|\zeta(1-1/\log\log x+it)|^2}{t\sqrt{(2\alpha t)^2+1}}\diff t\ll \int_{(\log x)^2}^\infty \frac{t^{1/10}}{t^{1+1/2}}\diff t\ll 1.
\end{align*}
For the remaining $2\leq t\leq (\log x)^2$, we instead use \eqref{eq:zeta_line} to prove that 
\begin{align*}
\int_2^{(\log x)^2}\frac{|\zeta(1-1/\log\log x+it)|^2}{t\sqrt{(2\alpha t)^2+1}}\diff t\ll \int_2^{(\log x)^2}\frac{(\log t)^2}{t}\diff t\ll (\log\log x)^3.
\end{align*}
Collecting the inequalities above, we find the expected
\begin{align*}
\numberthis\label{eq:maj_Cal_I}\mathcal I(x,m,\alpha)&\ll \sqrt x \frac{m^{1/\log(1+m)^{1/4}}}{m_0}(\log\log x)^3.
\end{align*}

Now if $\alpha_0\leq \alpha<\alpha_1\ll 1/\log x$, then $1\leq X_\alpha<2$. We can then use \eqref{mainterm2} to deduce that 
\begin{align*}
\sum_{d\in\D_\alpha(x)}\chi_d(m)&=(\sqrt x-\sqrt{Y_1(1,\alpha)})\frac{(-1)^{\omega(m_0)}}{m_0}\prod_{\substack{2\leq j\leq r\\ e_j~even}}\left(1-\frac2{p_j}\right)\frac{a(m)}4+\O\left(m x^{\alpha}\log x\right)
\\&=(\sqrt x-\sqrt{Y_1(1,\alpha)})\EE(\XX_\infty(m))+\O\left(m\log x\right).
\end{align*}
The second equality relies on \nameref{lemm:esp}. This proves the last part of our proposition.
\end{proo}

\section{Mean value of \texorpdfstring{$L(1,\chi_d)^z$}{Dirichlet L-functions}}
We need a good approximation for $L(1,\chi_d)^z$, and the one given by Dahl-Lamzouri appears to be perfectly suitable.
\begin{prop}[Proposition 3.3 of \cite{dahl_distribution_2016}]\label{prop:appro-L} Let $d\in\D$ be large and let $0<\varepsilon<1/2$ be fixed. Let $y$ be a real number such that $\log d/\log_2d\leq \log y\leq \log d$. Assume that $L(s,\chi_d)$ has no zeros inside the rectangle $\{s\in\CC:1-\varepsilon<\Re(s)\leq 1, |\Im(s)|\leq 2(\log d)^{2/\varepsilon}\}$. Then for any complex $z$ such that $|z|\leq \log y/(4\log_2d\log_3d)$ we have that
\[L(1,\chi_d)^z=\sum_{n=1}^{\infty}\frac{ d_z(n)\chi_d(n)}ne^{-n/y}+\O_\varepsilon\left(\exp\left(-\frac{\log y}{2\log_2d}\right)\right).\]
Here, we recall that $d_z(n)$ are defined to be the coefficients of the Dirichlet series of $\zeta^z$.
\end{prop}
We will use this estimate to prove our \nameref{prop:sum_Lz}. First we shall mention the following standard bound (see Lemma 2.2 of \cite{lamzouri_extreme_2010}): if $\chi$ is a non-exceptional and non-principal Dirichlet character modulo $q$, then for all $t\in\RR$
\[\numberthis\label{eq:bound-L}\log L(1+it,\chi)\ll\log_2(q(|t|+2)).\]

\begin{proo}[of \nameref{prop:sum_Lz}]
Fix $\frac{4(\alpha+1)}{2\alpha+5}<\sigma_\alpha<1$, say $\sigma_\alpha:=5/6+\alpha/3$. Let $M$ be the number of positive fundamental discriminants $d\leq x$ such that $L(s,\chi_d)$ has a zero in the rectangle
\[\R(x):=\{s\in\CC: \sigma_\alpha<\Re(s)\leq 1,|\Im(s)|\leq 2(\log x)^{2/(1-\sigma_\alpha)}\}.\]
Then, we use Heath-Brown's zero-density estimate (see \cite{heath-brown_mean_1995}), which states that for any $\varepsilon>0$ and $1/2<\sigma<1$, we have that 
\[\sum_{|d|\leq x}{\vphantom{\sum}}^{\!\!\flat} N(\sigma,T,\chi_d)\ll(xT)^\varepsilon x^{3(1-\sigma)/(2-\sigma)}T^{(3-2\sigma)/(2-\sigma)},\]
where $N(\sigma,T,\chi_d)$ is the number of zeros $\rho$ of $L(s,\chi_d)$ with $\Re(\rho)\geq \sigma$ and $|\Im(\rho)|\leq T$, and $\sum^\flat$ indicates that the sum is over fundamental discriminants. Thus, for any small $\varepsilon>0$ we have that 
\[\numberthis\label{majM} M\ll x^{3(1-\sigma_\alpha)/(2-\sigma_\alpha)+\varepsilon}.\]

Now we have to link this number $M$ to the number $N$ of $d\in\D_\alpha(x)$ for which $L(s,\chi_d)$ has a zero in $\R(x)$. We let $d_1,...,d_M$ be the positive fundamental discriminants lower than $x$ for which $L(s,\chi_d)$ has a zero in $\R(x)$. We recall that for every $d\in\D$, there is a unique positive fundamental discriminant $\tilde d$ and some $\ell\in\NN$ such that $d=\tilde d\ell^2$. In this case, it is easy to see that $L(s,\chi_d)$ and $L(s,\chi_{\tilde d})$ have the same zeros in the region $\Re(s)>0$.

Consequently, if we fix $d\in\D_\alpha(x)$ such that $L(s,\chi_d)$ has a zero in $\R(x)$, we may write that for some $1\leq u\leq X_\alpha,\,  Y_2\leq t\leq Y_3, \,j\in\{1,...,M\},\, \ell\in\NN$, we have $\frac{t^2-4}{u^2}=d_j\ell^2$, \ie $t^2-(u^2d_j)\ell^2=4$. Therefore, by the theory of Pell's equation there exists $n\in\NN$ such that
\[\frac{t+\ell\sqrt{d_j}u}2=\varepsilon_{d_ju^2}^n.\]
Thus, for any fixed $1\leq u\leq X_\alpha$, the fact that $\varepsilon_{d_ju^2}^n\leq t\leq Y_3\ll x$ yields
\[|\{(t,\ell): Y_2\leq t\leq Y_3, \ell\in\NN: t^2-4=u^2\ell^2d_j\}|\leq|\{n\in\NN:\varepsilon_{d_ju^2}^n\ll x\}|\ll\log x.\]
By \eqref{majM} it follows that 
\[\numberthis\label{eq:majM}N \ll M (\log x) X_\alpha\ll x^{3(1-\sigma_\alpha)/(2-\sigma_\alpha)+\alpha+\varepsilon},\]
for any small $\varepsilon>0$. Note that 
\[\numberthis\label{eq:majM2}3\frac{1-\sigma_\alpha}{2-\sigma_\alpha}+\alpha=\frac12-\frac{(1-2\alpha)^2}{14-4\alpha}.\]
We define $\delta_\alpha:=(1-2\alpha)^2/15$ so that \eqref{eq:majM} and \eqref{eq:majM2} imply that
\[ N\ll x^{1/2-\delta_\alpha}.\]

Let $\tilde{D}_\alpha(x)$ be the set of $d\in\D_\alpha(x)$ such that $d> \sqrt x$ and $L(s,\chi_d)$ has no zeros in $\R(x)$, then the last inequality yields
\[\numberthis\label{eq:dif-d}|\D_\alpha(x)|-|\tilde{\D}_\alpha(x)|\ll x^{1/2-\delta_\alpha}.\]
Thus, using \eqref{eq:bound-L}, we may write 
\[\sum_{d\in\D_\alpha(x)}L(1,\chi_d)^z-\sum_{d\in\tilde\D_\alpha(x)}L(1,\chi_d)^z\ll x^{1/2-\delta_\alpha}\exp(\O(|z|\log_2 x))\ll x^{1/2-\delta_\alpha/2}.\]
Let $y:=x^{\delta_\alpha}$ and $k:=\lceil |z|\rceil$, so that \nameref{prop:appro-L} and \eqref{eq:card_D} imply that
\begin{align*}
&\sum_{d\in\D_\alpha(x)}L(1,\chi_d)^z=\sum_{d\in\tilde\D_\alpha(x)}L(1,\chi_d)^z+\O(x^{1/2-\delta_\alpha/2})
\\&=\sum_{d\in\tilde{\D}_\alpha(x)}\sum_{m=1}^\infty \frac{d_z(m)\chi_d(m)e^{-m/y}}m+\O\left(|\tilde D_\alpha(x)|\exp\left(\frac{-\delta_\alpha\log x}{2\log_2 x}\right)\right)+\O(x^{1/2-\delta_\alpha/2})
\\&=\numberthis\label{exp:sum-L}\sum_{d\in\tilde{\D}_\alpha(x)}\sum_{m=1}^\infty \frac{d_z(m)\chi_d(m)e^{-m/y}}m+\O\left(\sqrt x\exp\left(\frac{-\delta_\alpha\log x}{3\log_2 x}\right)\right).
\end{align*}
We now extend the main term to a sum over all $d\in\D_\alpha(x)$. We may do so, for \eqref{eq:dz} and \eqref{eq:dif-d} imply that
\[\sum_{d\in\D_\alpha(x)\backslash\tilde\D_\alpha(x)}\sum_{m=1}^\infty\frac{d_z(m)\chi_d(m)}me^{-m/y}\ll (|\D_\alpha(x)|-|\tilde{\D}_\alpha(x)|)\sum_{m=1}^\infty\frac{d_k(m)e^{-m/y}}{m}\ll x^{1/2-\delta_\alpha/2}.\]
Combining this with \eqref{exp:sum-L}, we get that
\[\numberthis\label{sum_estim_L}\sum_{d\in\D_\alpha(x)}L(1,\chi_d)^z=\sum_{m=1}^\infty \frac{d_z(m)e^{-m/y}}m\sum_{d\in\D_\alpha(x)}\chi_d(m)+\O\left(\sqrt x\exp\left(\frac{-\delta_\alpha\log x}{3\log_2 x}\right)\right).\]

We first assume that $\alpha_1\leq \alpha\leq \alpha'$. We may use \nameref{prop:sum-chid}, denoting its main term by 
\begin{align*}
&MT(m):=MT(x,m):=
\\&\sqrt x\Biggl(\frac{\alpha^2}{2\zeta(2)}\EE(\XX(m))\log^2x+\alpha\left[H'_m(1)+(2\gamma-2\alpha)\frac{\EE(\XX(m))}{\zeta(2)}\right](\log x-2)+\eth_2\Biggr),
\end{align*}
to write that
\begin{align*}
\sum_{m=1}^\infty \frac{d_z(m)e^{-m/y}}m&\sum_{d\in\D_\alpha(x)}\chi_d(m)=\sum_{m=1}^\infty \frac{d_z(m)e^{-m/y}}{m}MT(m)
\\&+(1-\sqrt{Y_1(1,\alpha)})\sum_{m=1}^\infty\frac{d_z(m)e^{-m/y}}m\EE(\XX_\infty(m))+\O\left(x^\alpha\log x\sum_{m=1}^\infty d_k(m)e^{-m/y}\right)
\\&\numberthis\label{sum:dz-chi}+\sum_{m=1}^\infty \frac{d_z(m)e^{-m/y}}{m}\left [\mathcal I(x,m,\alpha)+ \O\left(\frac{\log^2(m+1)}{m_0}\right)\right].
\end{align*}

The aim is now to get rid of the weigh $e^{-m/y}$ in the main term. By \nameref{lemm:esp} and \nameref{majoration_MT}, we know that removing the said weigh introduces an error bounded by
\[\sqrt x\log^2x\sum_{m=1}^\infty \frac{d_k(m)(1-e^{-m/y})\log^2(m+1)}{m_0m}.\]
We know that $1-e^{-t}\ll t^{\beta}$ for all $t>0$ and $0<\beta<1$, and we choose $\beta=1/\log_2 x$. By \nameref{lemm:majoration-m0}, if $|z|\leq\frac{\delta_\alpha}5\frac{\log x}{\log_2 x\log_3 x}$, we deduce that this error is bounded by 
\begin{align*}
\sqrt xy^{-\beta}\log^2x \sum_{m=1}^\infty \frac{d_k(m)m^\beta\log^2(m+1)}{m_0m}&\ll \sqrt x x^{-\delta_\alpha/\log_2 x} \exp\left(\frac{\delta_\alpha}2\frac{\log x}{\log_2 x}\right)
\\&\numberthis\label{eq:remove_weight}\ll \sqrt{x}\exp\left(-\frac{\delta_\alpha}2\frac{\log x}{\log_2 x}\right).
\end{align*}
Similarly since $\mathcal{I}(x,m,\alpha)\ll \sqrt x \frac{m^{1/\log_2(2+ m)^{1/4}}}{m_0}(\log\log x)^3$ by \nameref{prop:sum-chid} and since $\EE(\XX_\infty(m))\ll 1/m_0$ by \nameref{lemm:esp} the error introduced by removing the exponential weight in the second and fourth sums of \eqref{sum:dz-chi} above is 
\[\numberthis\label{eq:I_final1}\ll \sqrt x\exp\left(-\frac{\delta_\alpha}2\frac{\log x}{\log_2 x}\right).\]
Now that the exponential weight is removed, \nameref{lemm:majoration-m0} shows that the contribution of the $\log(1+m)^2/m_0$ term in \eqref{sum:dz-chi} is 
\[\numberthis\label{eq:I_final2}\ll \exp\left(\O\left(\frac{\log x}{\log_2 x}\right)\right).\]
Since $Y_1(1,\alpha)\ll x$ and by \nameref{majoration_esp_s}, the contribution of the $\EE(\XX_\infty(m))$ term in \eqref{sum:dz-chi} is 
\[(1-\sqrt{Y_1(1,\alpha)})\EE(L(1,\XX_\infty)^z)\ll \sqrt x \EE(L(1,\XX)^{\Re(z)}).\]
Finally, \nameref{lemm:majoration_sans_m} ensures that the contribution of the error term of the third sum in \eqref{sum:dz-chi} is 
\[\numberthis\label{eq:ger-err-term}\ll x^\alpha(\log x) y(\log(3y))^{k+2}\ll x^{\alpha+2\delta_\alpha}\ll x^{1/2-\delta_\alpha/2}.\]

Collecting all these estimates and using \eqref{sum_estim_L}, we find that 
\begin{align*}
\sum_{d\in\D_\alpha(x)}L(1,\chi_d)^z&=\sum_{m=1}^\infty \frac{d_z(m)}{m}MT(m)+\sum_{m=1}^\infty \frac{d_z(m)}{m}\mathcal I(x,m,\alpha)
\\\numberthis\label{eq:aux_stiel}&\qquad+(1-\sqrt{Y_1(1,\alpha)})\EE(L(1,\XX_\infty)^z)+\O\left(\exp\left(-\frac{-\delta_\alpha}3\frac{\log x}{\log_2 x}\right)\right)
\\\numberthis\label{eq:expression_finale}&=\sum_{m=1}^\infty \frac{d_z(m)}{m}MT(m)+\sum_{m=1}^\infty \frac{d_z(m)}{m}\mathcal I(x,m,\alpha)+\O\left(\sqrt x \EE(L(1,\XX)^{\Re(z)})\right).
\end{align*}
Using \nameref{lemm:exp_phi}, we find that the sum of $MT(m)$ is equal to
\begin{align*}
&\sum_{m=1}^\infty \frac{d_z(m)MT(m)}{m}
\\&=\EE(L(1,\XX)^z)\frac{\sqrt x}{\zeta(2)}\Biggl(\alpha^2\frac{\log^2x}{2}+\alpha\left(\phi(z)-2\alpha\right)(\log x-2)\Biggr)+\sqrt x\sum_{m=1}^\infty\frac{d_z(m)}m \eth_2(m)
\\&=\numberthis\label{eq:sum_mt}\EE(L(1,\XX)^z)\frac{\sqrt x}{\zeta(2)}\Biggl(\alpha^2\frac{\log^2x}{2}+\alpha\left(\phi(z)-2\alpha\right)(\log x-2)\Biggr)+\O\left(\sqrt x\EE(L(1,\XX)^{\Re(z)}\right).
\end{align*}
The last equality relies on \nameref{coro:majo_esp}, and on the definition \eqref{eq:defEth} of $\eth_2$.

It only remains to estimate the contribution of $\mathcal I(x,m,\alpha)$. Observe that, for $s=\sigma+i\tau\in\CC$, $3/2<\sigma<2$ : $1/|\zeta(s)|<\zeta(\sigma)/\zeta(2\sigma)\ll 1$, and hence on the line $\Re(s)= -1/\log\log x$, we have $\frac{1}{\zeta(2(s+1))}\ll1$.
Furthermore, as seen through the proof of \nameref{lemm:exp_phi}, we have $\sum_{m=1}^\infty\frac{d_z(m)}mH_m(s)=\frac{\psi(s)}{\zeta(2s)}\EE(L(1,\XX_s)^z)$, where $\psi(s)\ll 1$ on the line $\Re(s)= 1-1/\log\log x$. By \nameref{majoration_esp_s} and using the same computations that lead to \eqref{eq:maj_Cal_I}, we deduce that 
\begin{align*}
\sum_{m=1}^\infty \frac{d_z(m)}{m}\mathcal I(x,m,\alpha)&=\frac{X_\alpha^{1/(2\alpha)}}{2\pi i}\int_{-1/\log\log x-i\infty}^{-1/\log\log x+i\infty}\frac{\psi(1+s)\EE(L(1,\XX_{1+s})^z)\zeta^2(1+s)}{\zeta(2(s+1))s(2\alpha s+1)}X_\alpha^s\diff s
\\&\numberthis\label{eq:inte_a_ref}\ll \sqrt x\EE(L(1,\XX)^{\Re(z)})(\log\log x)^3.
\end{align*}
Using the estimate above and \eqref{eq:sum_mt} in \eqref{eq:expression_finale}, we have proved that if $\alpha_1\leq \alpha\leq \alpha'$, then
\begin{align*}
\sum_{d\in\D_\alpha(x)}L(1,\chi_d)^z&=\EE(L(1,\XX)^z)\frac{\sqrt x}{\zeta(2)}\Biggl(\alpha^2\frac{\log^2x}{2}+\alpha\left(\phi(z)-2\alpha\right)(\log x-2)\Biggr)
\\&\qquad\numberthis\label{eq:1Lpetit}+\O\left(\sqrt x(\log\log x)^3\EE(L(1,\XX)^{\Re(z)}\right).
\end{align*}

Now we assume that $\alpha_0\leq \alpha<\alpha_1$. Then we can go back to \eqref{sum_estim_L} and use \nameref{prop:sum-chid} to find that 
\begin{align*}
\sum_{d\in\D_\alpha(x)}L(1,\chi_d)^z&=\sum_{m=1}^\infty \frac{d_z(m)e^{-m/y}}m\sum_{d\in\D_\alpha(x)}\chi_d(m)+\O\left(\sqrt x\exp\left(\frac{-\delta_\alpha\log x}{3\log_2 x}\right)\right)
\\&=(\sqrt x-\sqrt{Y_1(1,\alpha)})\sum_{m=1}^\infty\frac{d_z(m)}{m}\EE(\XX_\infty(m))+\O\left(\sum_{m=1}^\infty d_k(m)e^{-m/y} \log x\right)
\\&\numberthis\label{eq:equation_petit_u}\qquad+\O\left(\sqrt x\exp\left(\frac{-\delta_\alpha\log x}{3\log_2 x}\right)\right)
\\&=\O(\sqrt x\EE(L(1,\XX)^{\Re(z)})),
\end{align*}
by \nameref{majoration_esp_s}. The exponential weight of the first term of the second line was removed as in \eqref{eq:remove_weight}. The contribution of the first error term of the second line was bounded as in \eqref{eq:ger-err-term}.

Combining the two results above, our \nameref{rema:phi} and the fact that $\alpha_1\ll 1/\log x$ show that \eqref{eq:1Lpetit} holds if $0<\alpha\leq \alpha_1$, which concludes the proof.
\end{proo}

\section{Proof of \nameref{theo-hd}}\label{sec:stiel}
We can now prove \nameref{theo-hd}, simply by using partial summation and Stieltjes integration. The fact that we want uniformity in the moments makes this step of the proof much more subtle compared to that of Hooley's paper \cite{hooley_pellian_1984} (see (63) page 123). We will work with a fixed $0<\alpha<1/2$ and a large real number $x$. We recall that $X_\alpha=x^\alpha-x^{-1-\alpha}$ and that $\alpha_0$ (resp. $\alpha_1$) is such that $X_{\alpha_0}=1$ (resp. $X_{\alpha_1}=2$). We fix a complex number $z$ such that $\sigma:=\Re(z)\geq -1$ and $|z|\leq\frac{(1-2\alpha)^2}{75}\frac{\log x}{\log_2x\log_3x}$, and we define $f(t):=(\sqrt t/\log t)^z$. Observe that since $|f'(t)|\ll |z| \left(\frac{\sqrt t}{\log t}\right)^{\sigma-1}\frac{1}{\sqrt t \log t}$, we may deduce that if $\sigma>-1$:
\begin{align*}
\numberthis\label{eq:f'1}\int_2^x \sqrt tf'(t)\diff t\ll |z|\int_2^x\frac{(\sqrt t)^{\sigma-1}}{(\log t)^\sigma}\diff t\ll \frac{|z|}{\sigma+1}\frac{x^{(\sigma+1)/2}}{(\log x)^\sigma}.
\end{align*}
If $\sigma=-1$, then 
\[\numberthis\label{eq:f'2}\int_2^x \sqrt tf'(t)\diff t\ll |z| \log^2 x.\]
\begin{proo}[of \nameref{theo-hd}]
First, we want to estimate
\[Q(x,u,z):=\sum_{d\in\D_u(x)}\left(\frac{\sqrt dL(1,\chi_d)}{\log d}\right)^z,\]
uniformly in $\alpha_0 \leq u\leq\alpha$. We first let $\alpha_0\leq u< \alpha_1$ and $2\leq t\leq x$. Going back to \eqref{eq:equation_petit_u}, we let 
\[MT_1(t,u,z):=(\sqrt t-Y_1(1,u))\EE(L(1,\XX_\infty)^z)\]
so that, proceeding as in \eqref{eq:ger-err-term} to bound the error term, we have 
\[\sum_{d\in\D_u(t)}L(1,\chi_d)^z=MT_1(t,u,z)+Err_1(t,u,z)\]
with $Err_1(t,u,z)\ll \sqrt t\exp(-\delta\log t/\log_2(t+1))$ for some $\delta>0$. We also put $MT_1(y,u,z)=0$ if $0\leq y<2$.

By partial summation and \nameref{majoration_esp_s}, this leads to
\begin{align*}
Q(x,u,z)&=\int_{2^-}^xf(t)\frac{\partial}{\partial t}MT_1(t,u,z)\diff t+Err_1(x,u,z)f(x)-\int_{2^-}^x Err_1(t,u,z)f'(t)\diff t
\\\numberthis\label{eq:Q_alph}&=\EE(L(1,\XX_\infty)^z)\int_2^x\frac{f(t)}{2\sqrt t}\diff t+Err_{Q,1}(x,u,z).
\end{align*}
By \eqref{eq:f'1} if $\sigma>-1$:
\[\numberthis\label{eq:maj_err_q1} Err_{Q,1}(x,u,z)\ll \left(\frac{|z|}{\sigma+1}+1\right)\frac{x^{(\sigma+1)/2}}{(\log x)^{\sigma+1}},\]
and by \eqref{eq:f'2}, if $\sigma=-1$:
\[\numberthis\label{eq:maj_err_q11}Err_{Q,1}(x,u,z)\ll |z|.\]

Now we let $\alpha_1\leq u \leq \alpha$ and $2\leq t\leq x$. We can not use \nameref{prop:sum_Lz} directly, and we instead need to go back to \eqref{eq:aux_stiel} and use \eqref{eq:sum_mt} to find that 
\begin{align*}
\sum_{d\in\D_u(t)}L(1,\chi_d)^z&=\EE(L(1,\XX)^z)\frac{\sqrt t}{\zeta(2)}\Biggl(u^2\frac{\log^2t}{2}+u\left(\phi(z)-2u\right)(\log t-2)\Biggr)+\sqrt t\sum_{m=1}^\infty\frac{d_z(m)}m\eth_2(m)
\\&+(1-\sqrt{Y_1(1,u)})\EE(L(1, \XX_\infty)^z)+\sum_{m=1}^\infty \frac{d_z(m)}m\mathcal I(t,m,u)
\\&\numberthis\label{eq:sum_L_precise}+\O\left(\EE(L(1,\XX)^\sigma)\sqrt t\exp\left(-\delta\frac{\log t}{\log_2(t+1)}\right)\right),
\end{align*}
for some $\delta>0$ and where $\eth_2(m)$ and $\mathcal I(t,m,u)$ were defined in \nameref{prop:sum-chid}. The reason for which we need this lengthier expression instead of the easier \nameref{prop:sum_Lz} is that the error term above is much more precise than the one in \nameref{prop:sum_Lz}. By keeping in the main term some expressions rather than putting them in the error term (as we did in \nameref{prop:sum_Lz}), we will be able to handle them with a lot more precision in the following partial summations, leading to a better final error term in \nameref{theo-hd}. However we can already show that some of the contribution of the term $\sum_{m=1}^\infty \frac{d_z(m)}m\mathcal I(t,m,u)$ goes in the error term of \eqref{eq:sum_L_precise}.

First, observe that
\[\numberthis\label{eq:derivee0}\frac{\partial }{\partial u}\frac{\partial}{\partial t}\left\{(1-\sqrt{Y_1(1,u)})\EE(L(1, \XX_\infty)^z)+\sqrt t\sum_{m=1}^\infty\frac{d_z(m)}m\eth_2(m)\right\}=0,\]
and we denote the term between curly brackets by $C(t,u,z)$. As proved at the very beginning of \nameref{lemm:exp_phi}, we know that $\sum_{m=1}^\infty\frac{d_z(m)}mH_m(s)=\frac{\psi(s)}{\zeta(2s)}\EE(L(1,\XX_s)^z)$, where $\psi(s)\ll 1$ if $\Re(s)$ is close to $1$. Therefore,
\[\sum_{m=1}^\infty \frac{d_z(m)}m\mathcal I(t,m,u)=\frac{(t^u-t^{-1-u})^{1/(2u)}}{2\pi i}\int_{(-c)}\frac{\psi(1+s)\EE(L(1,\XX_{1+s})^z)\zeta^2(1+s)}{\zeta(2(s+1))s(2us+1)}(t^u-t^{-1-u})^s\diff s.\]
for any small enough $c>0$, where $(-c)$ means that the integration path is the line $\Re(s)=-c$, and we choose $c=1/\log\log x$. For commodity, we put $\Psi_z(s):=\frac{\psi(s)\EE(L(1,\XX_s)^z)}{2\pi i\zeta(2s)}$, and hence on the line $\Re(s)=-c$
\[\numberthis\label{eq:maj_Psi}\Psi_z(1+s)\ll \EE(L(1,\XX)^z)\]
by \nameref{majoration_esp_s} and since $1/\zeta(2(s+1))\ll 1$, as shown earlier, above \eqref{eq:inte_a_ref}.

Since $\alpha_1\leq u\leq \alpha$, we may use \eqref{eq:devAsympX} and bound the integral as we did in \eqref{eq:inte_a_ref} to write that 
\[\sum_{m=1}^\infty \frac{d_z(m)}m\mathcal I(t,m,u)=\sqrt t\int_{(-c)}\frac{\Psi_z(1+s)}{s(2su+1)}\zeta^2(1+s)(t^u-t^{-1-u})^s\diff s+\O\left(\EE(L(1,\XX)^\sigma)t^{-1/3}\log\log^3 x \right) .\]
We fix $\varepsilon>0$. Since $u\gg 1/\log x$, we have by \eqref{eq:maj_Psi} and the Phragmén-Lindelöf theorem \eqref{eq:mu}
\[\sqrt t\int_{\substack{(-c)\\|\Im(s)|>\log^3x}}\frac{\Psi_z(1+s)}{s(2su+1)}\zeta^2(1+s)(t^u-t^{-1-u})^s\diff s\ll \sqrt t\EE(L(1,\XX)^\sigma)\int_{\log^3 x}^\infty \frac{t^\varepsilon}{t^2u}\diff t\ll \sqrt t\frac{\EE(L(1,\XX)^\sigma)}{\log x}.\]
By \eqref{eq:logXa}, for $s\in\CC$ with $\Re(s)=-c$ and $|\Im(s)|\leq \log^3x$, we have $(t^u-t^{-1-u})^s=t^{us}+\O(1/\sqrt t)$. Thus by using \eqref{eq:maj_Psi}, the fact that $\zeta(1-c+iw)\ll \log w$ if $2\leq w\leq \log x^3$ by \eqref{eq:zeta_line} and the fact that $\zeta(1-c+it)\ll 1/c$ if $|t|\leq 2$:
\begin{align*}
\sum_{m=1}^\infty \frac{d_z(m)}m\mathcal I(t,m,u)&=\sqrt t\int_{\substack{(-c)\\|\Im(s)|\leq \log^3x}}\frac{\Psi_z(1+s)\zeta^2(1+s)}{s(2us+1)}t^{us}\diff s
\\&\qquad+\O\left(\EE(L(1,\XX)^\sigma)\left[\frac{1}{c^3}+\int_{\substack{(-c)\\2\leq |\Im(s)|\leq \log^3x}}\frac{\log^2w}{w}\diff w+\frac{\sqrt t}{\log x}\right]\right)
\\=\sqrt t\int_{\substack{(-c)\\|\Im(s)|\leq \log^3x}}&\frac{\Psi_z(1+s)\zeta^2(1+s)}{s(2us+1)}t^{us}\diff s+\O\left(\EE(L(1,\XX)^\sigma)\left[\log\log^3 x+\frac{\sqrt t}{\log x}\right]\right).
\end{align*}

Thus, going back to \eqref{eq:sum_L_precise}, we may define 
\begin{align*}
MT_2(t,u,z)&:=\EE(L(1,\XX)^z)\frac{\sqrt t}{\zeta(2)}\Biggl(u^2\frac{\log^2t}{2}+u\left(\phi(z)-2u\right)(\log t-2)\Biggr) +C(t,u,z)
\\&\qquad+\sqrt t\int_{\substack{(-c)\\|\Im(s)|\leq \log^3x}}\frac{\Psi_z(1+s)\zeta^2(1+s)}{s(2us+1)}t^{us}\diff s
\end{align*}
and $Err_2(t,u,z)$ by 
\[\sum_{d\in\D_u(t)}L(1,\chi_d)^z=MT_2(t,u,z)+Err_2(t,u,z),\]
so that the crude bound $\sqrt t\exp(-\delta\log t/\log\log(t+1))\ll t/\log t$ leads to
\[Err_2(t,u,z)\ll \EE(L(1,\XX)^\sigma)\left(\frac{\sqrt t}{\log t}+\log\log^3 x\right).\]

Doing the partial summation as we did in \eqref{eq:Q_alph}, we find that
\begin{align*}
&Q(x,u,z)
\\&=\frac{\EE(L(1,\XX)^z)}{\zeta(2)}\int_2^xf(t)\Biggl(\frac{u^2 \frac{\log^2t}2+u\left(\phi(z)-2u\right)(\log t-2)}{2\sqrt t}+\frac{u^2\log t+u\left(\phi(z)-2u\right)}{\sqrt t}\Biggr)\diff t
\\&\qquad+\int_2^x f(t)\frac{\partial}{\partial t}C(t,u,z)\diff t+\int_2^x \frac{f(t)}{2\sqrt t}\int_{\substack{(-c)\\|\Im(s)|\leq \log^3x}}\frac{\Psi_z(1+s)\zeta^2(1+s)}{s}t^{us}\diff s+ Err_{Q,2}(x,u,z)
\\&= \frac{\EE(L(1,\XX)^z)}{\zeta(2)}\int_2^x\frac{f(t)}{\sqrt t}\left\{\frac{u^2}4\log^2t+\frac{u}2\phi(z)\log t\right\}\diff t +\int_2^x f(t)\frac{\partial}{\partial t}C(t,u,z)\diff t
\\&\qquad\numberthis\label{eq:exp_Q_u_grand}+\int_2^x \frac{f(t)}{2\sqrt t}\int_{\substack{(-c)\\|\Im(s)|\leq \log^3x}}\frac{\Psi_z(1+s)\zeta^2(1+s)}{s}t^{us}\diff s+Err_{Q,2}(x,u,z),
\end{align*}
so that by \eqref{eq:f'1}
\[\numberthis\label{eq:err_1:stiel} Err_{Q,2}(x,u,z)\ll \left(\frac{|z|}{\sigma+1}+1\right)\frac{x^{(1+\sigma)/2}}{(\log x)^{\sigma+1}}\EE(L(1,\XX)^\sigma)\]
if $\sigma>-1$, and if $\sigma=-1$ \eqref{eq:f'2} implies that
\[\numberthis\label{eq:err_2:stiel}Err_{Q,2}(x,u,z)\ll |z|\log x\]
We define $MT_{Q,i}(x,u,z):=Q(x,u,z)-Err_{Q,i}(x,u,z)$, $i=1,2$, and $MT_{Q,i}(x,v,z)=0$ if $0\leq v<\alpha_0$, for convenience. 

Thanks to Dirichlet's Class Number Formula, we may write that
\begin{align*}
S(x,\alpha,z):=\sum_{d\in\D_\alpha(x)}h(d)^z&=\sum_{\substack{d\in\D(x)\\\varepsilon_d\leq d^{1/2+\alpha}}}\left(\frac{\sqrt dL(1,\chi_d)}{\log d}\right)^z\left(\frac{\log d}{\log\varepsilon_d}\right)^z
\\&=\int_{0}^\alpha \left(1/2+u\right)^{-z}\diff Q(x,u,z).
\end{align*}
We first study the integral over $(0,\alpha_1)$. By integrating by parts and using \eqref{eq:Q_alph} together with \eqref{eq:maj_err_q1} and \eqref{eq:maj_err_q11}, we find that
\begin{align*}
\int_{0}^{\alpha_1^-} \left(1/2+u\right)^{-z}\diff Q(x,u,z)&=\Bigl(MT_{Q,1}(x,\alpha_1,z)+Err_{Q,1}(x,\alpha_1,z)\Bigr)(1/2+\alpha_1)^{-z}
\\&\qquad-\int_{0}^{\alpha_1^-} Q(x,u,z)\frac{\partial}{\partial u}\left(1/2+u\right)^{-z}\diff u
\\&=\int_{0}^{\alpha_1^-} \left(1/2+u\right)^{-z}\frac{\partial}{\partial u}MT_{Q,1}(x,u,z)\diff u+Err_{S,1}(x,\alpha,z)
\\&\numberthis\label{equ:Q_stiel_petit}=Err_{S,1}(x,\alpha,z),
\end{align*}
where 
\begin{align*}
Err_{S,1}(x,\alpha,z) &= Err_{Q,1}(x,\alpha_1,z)\left(1/2+\alpha_1\right)^{-z}-\int_{0}^{\alpha_1^-} Err_{Q,1}(x,u,z)\frac{\partial}{\partial u}\left(1/2+u\right)^{-z}\diff u
\\&\ll 2^\sigma\left(\frac{|z|}{\sigma+1}+1\right)^2\frac{x^{(\sigma+1)/2}}{(\log x)^{\sigma+1}}
\end{align*}
if $\sigma>-1$, and if $\sigma=-1$ 
\[Err_{S,1}(x,\alpha,z)\ll |z|^2.\]

Now we deal with the integral over $[\alpha_1,\alpha)$. By similar computations and using \eqref{eq:derivee0}, \eqref{eq:exp_Q_u_grand}
\begin{align*}
&\int_{\alpha_1^+}^{\alpha} \left(1/2+u\right)^{-z}\diff Q(x,u,z)=\Bigl[Q(x,u,z)(1/2+u)^{-z}\Bigr]_{\alpha_1^+}^\alpha-\int_{\alpha_1^+}^{\alpha} Q(x,u,z)\frac{\partial}{\partial u}\left(1/2+u\right)^{-z}\diff u
\\&=\int_{\alpha_1}^\alpha \left(1/2+u\right)^{-z}\frac{\partial}{\partial u}MT_Q(x,u,z)\diff u+Err_{S,2}(x,\alpha,z)
\\&=\frac{\EE(L(1,\XX)^z)}{2\zeta(2)}\int_2^x\frac{t^{\frac{z-1}{2}}}{(\log t)^z}\left\{\left(\int_{\alpha_1}^\alpha u\left(1/2+u\right)^{-z}\diff u\right)\log^2t+\left(\int_{\alpha_1}^\alpha \left(1/2+u\right)^{-z}\diff u\right)\phi(z)\log t\right\}\diff t
\\&\qquad\numberthis\label{equ:inte_stiel_q}+\int_2^x \frac{f(t)\log t}{2\sqrt t}\int_{\substack{(-c)\\|\Im(s)|\leq \log^3x}}\Psi_z(1+s)\zeta^2(1+s)\int_{\alpha_1}^{\alpha} \left(1/2+u\right)^{-z}t^{us}\diff u\diff s+Err_{S,2}(x,\alpha,z),
\end{align*}
where 
\begin{align*}
Err_{S,2}(x,\alpha,z) &= [Err_{Q,2}(x,u,z)\left(1/2+u\right)^{-z}]_{\alpha_1^+}^\alpha-\int_{0}^\alpha Err_{Q,2}(x,u,z)\frac{\partial}{\partial u}\left(1/2+u\right)^{-z}\diff u.
\end{align*}
By \eqref{eq:err_1:stiel}, if $\sigma>-1$:
\[Err_{S,2}(x,\alpha,z) \ll 2^\sigma\EE(L(1,\XX)^\sigma)\left(\frac{|z|}{\sigma+1}+1\right)^2\frac{x^{(\sigma+1)/2}}{(\log x)^{\sigma+1}},\]
and by \eqref{eq:err_2:stiel}, if $\sigma=-1$:
\[Err_{S,2}(x,\alpha,z) \ll |z|^2\log x.\]

Now we want to extend the inner integrals to the whole interval $(0,\alpha)$, and we may do so for \nameref{rema:phi} together with the fact that $\alpha_1\ll 1/\log x$ imply that 
\begin{align*}
&\frac{\EE(L(1,\XX)^z)}{2\zeta(2)}\int_2^x\frac{t^{\frac{z-1}{2}}}{(\log t)^z}\left\{\left(\int_{0}^{\alpha_1} u\left(1/2+u\right)^{-z}\diff u\right)\log^2t+\left(\int_{0}^{\alpha_1} \left(1/2+u\right)^{-z}\diff u\right)\phi(z)\log t\right\}\diff t
\\&\numberthis\label{equ:tout_a}\ll 2^\sigma\EE(L(1,\XX)^\sigma)\int_{2}^x \frac{t^{(\sigma-1)/2}}{(\log t)^\sigma}\diff t\ll \frac{2^\sigma\EE(L(1,\XX)^\sigma)}{\sigma+1}\frac{x^{(\sigma+1)/2}}{(\log x)^\sigma}
\end{align*}
if $\sigma>-1$, and if $\sigma=-1$ the quantity above is 
\[\ll (\log x)^2.\]

Observe that for $2\leq t\leq x$ and $\Re(s)=-c$, $|\Im(s)|\leq \log^3x$, a summation by parts leads to:
\begin{align*}
\int_{\alpha_1}^{\alpha} \left(1/2+u\right)^{-z}t^{us}\diff u
&=\left[\frac{(1/2+u)^{-z}t^{us}}{s\log t}\right]_{\alpha_1}^\alpha+\frac{z}{s\log t}\int_{\alpha_1}^\alpha (1/2+u)^{-z-1}t^{us}\diff u
\\&\ll \frac{2^\sigma}{|s|\log t}\left(1+\frac{|z|}{\sigma+2}\right).
\end{align*}
Thus, once again using \eqref{eq:zeta_line}, \eqref{eq:maj_Psi} and the fact that $\zeta(1+s)\ll 1/|s|$ when $s\to 0$:
\begin{align*}
&\int_2^x \frac{f(t)\log t}{2\sqrt t}\int_{\substack{(-c)\\|\Im(s)|\leq \log^3x}}\Psi_z(1+s)\zeta^2(1+s)\int_{\alpha_1}^{\alpha} \left(1/2+u\right)^{-z}t^{us}\diff u\diff s\diff t
\\&\hspace{1.4cm}\ll 2^\sigma\EE(L(1,\XX)^\sigma)\left(1+\frac{|z|}{\sigma+2}\right)\int_2^x\frac{\sqrt t^{\sigma-1}}{(\log t)^\sigma}\int_{\substack{(-c)\\|\Im(s)|\leq \log^3x}}\frac{|\zeta^2(1+s)|}{|s|}|\diff s|\diff t
\\&\hspace{1.4cm} \ll 2^\sigma\EE(L(1,\XX)^\sigma)\left(1+\frac{|z|}{\sigma+2}\right)\left(\int_2^x \frac{\sqrt t^{\sigma-1}}{(\log t)^\sigma}\diff t\right)\left(\frac{1}{c^3}+\int_{\substack{(-c)\\2\leq |\Im(s)|\leq \log^3x}}\frac{|\zeta^2(1+s)|}{|s|}|\diff s|\right)
\\& \hspace{1.4cm} \ll\frac{2^\sigma\EE(L(1,\XX)^\sigma)}{\sigma+1}\frac{x^{(\sigma+1)/2}}{(\log x)^\sigma}\left(1+\frac{|z|}{\sigma+2}\right)\log\log^3 x
\end{align*}
if $\sigma>-1$, while if $\sigma=-1$ this quantity is
\[\ll|z|(\log^2x)\log\log^3x.\]

We use the definition of $I_{z, j}(\alpha)$, $j=0,1$, given in \nameref{theo-hd} so that using the result above together with \eqref{equ:Q_stiel_petit}, \eqref{equ:inte_stiel_q} and \eqref{equ:tout_a}, we find that
\begin{align*}
S(x,\alpha,z)=\frac{\EE(L(1,\XX)^z)}{2\zeta(2)}\int_2^x\frac{t^{\frac{z-1}{2}}}{(\log t)^z}\Bigl[I_{z,1}(\alpha)\log^2t+I_{z,0}(\alpha)\phi(z)\log t\Bigr]\diff t+Err
\end{align*}
where, if $\sigma>-1$:
\begin{align*}
Err&\ll 2^\sigma\EE(L(1,\XX)^\sigma)\frac{x^{(\sigma+1)/2}}{(\log x)^{\sigma}}\left(\frac{\left(\frac{|z|}{\sigma+1}+1\right)^2}{\log x}+\frac{\log\log^3x}{\sigma+1}\left(1+\frac{|z|}{\sigma+2}\right)\right)
\\&\ll 2^\sigma\EE(L(1,\XX)^\sigma)\frac{x^{(\sigma+1)/2}}{(\log x)^{\sigma}}\left(\frac{|z|^2}{(\sigma+1)^2\log x}+\frac{\log\log^3x}{\sigma+1}\left(1+\frac{|z|}{\sigma+2}\right)\right)
\end{align*}
and if $\sigma=-1$:
\[Err\ll |z|(\log^2 x)\log\log^3x.\]
\end{proo}

\begin{proo}[of \nameref{coro_equiv}]
We fix $z$ a complex number with a large real part $\Re(z)=\sigma$, and we fix $0<\alpha<1/2$. By \nameref{rema:phi}, we know that
\[\numberthis\label{rappelPhi}\EE(L(1,\XX)^z)\phi(z)\ll \EE(L(1,\XX)^\sigma).\]
Furthermore, if we assume that $\Im(z)\ll \sigma$ and $\sigma\to \infty$, then:
\begin{align*}
\frac{I_{z,0}(\alpha)}{I_{z,1}(\alpha)}&=\frac{2}{\frac{1-z}{2-z}\frac{(2\alpha+1)^{2-z}-1}{(2\alpha+1)^{1-z}-1}-1}=\frac{2}{\frac{z-1}{z-2}\left(2\alpha+1+\frac{2\alpha}{(2\alpha+1)^{1-z}-1}\right)-1}
\\&=\frac{2}{\frac{z-1}{z-2}\left(2\alpha+1-2\alpha+\O\left((2\alpha+1)^{-\sigma}\right)\right)-1}
\\&=\frac{2}{\frac{1}{z-2}+\O((2\alpha+1)^{-\sigma})}\ll \sigma\ll \log x/(\log_2 x\log_3 x).
\end{align*}
Using \nameref{theo-hd} with \eqref{rappelPhi}, this proves the first part of the Corollary:
\begin{align*}
\sum_{d\in\D_\alpha(x)}h(d)^z=\frac{1}{\zeta(2)}\frac{x^{\frac{z+1}2}}{z+1}I_{z,1}(\alpha)\log(x)^{-z+2}\left(\EE(L(1,\XX)^z)+\O\left(\frac{\EE(L(1,\XX)^\sigma)}{\log_2x\log_3 x}\right)\right)+Err
\end{align*}
Now if $z$ is large and real and satisfies $z\leq \frac{(1-2\alpha)^2}{75}\frac{\log x}{(\log_2x)^2}$, we have that $I_{z,1}(\alpha)\sim 2^z/(4z^2)$ which implies that 
\begin{align*}
&\sum_{d\in\D_\alpha(x)}h(d)^z=\frac{1}{\zeta(2)}\frac{\EE(L(1,\XX)^z)}{z+1}I_{z,1}(\alpha)x^{\frac{z+1}2}\log(x)^{-z+2}\left(1+\O\left(\frac1{\log_2x\log_3 x}\right)\right)+Err
\\&=\frac{1}{\zeta(2)}\frac{\EE(L(1,\XX)^z)}{z+1}I_{z,1}(\alpha)x^{\frac{z+1}2}\log(x)^{-z+2}\left(1+\O\left(\frac1{\log_2 x\log_3 x}\right)+\O\left(\frac{z^2}{\log^2 x}(\log\log x)^3\right)\right),
\end{align*}
which proves the result.
\end{proo}

\section{Distribution of class numbers over \texorpdfstring{$\D_\alpha(x)$}{Hooley's family}.}
To prove \nameref{theo:distribution}, we will require an estimate of $\EE(L(1,\XX)^r)$ when $r$ is a large real number.
\begin{prop}[Proposition 4.2 of \cite{dahl_distribution_2016}]\label{prop:size_H}
For any real number $r\geq 4$, we have that
\[\log\EE(L(1,\XX)^r)=r\left(\log_2 r+\gamma+\frac{C_0-1}{\log r}+\O\left(\frac{1}{(\log r)^2}\right)\right).\]
\end{prop}
Once again, the proof found in \cite{dahl_distribution_2016} remains true for our random model, since $a_p,b_p=1/2+\O(1/p)$ and $a_p-b_p,c_p\ll1/p$. We are now ready to prove \nameref{theo:distribution}.
\begin{proo}[of \nameref{theo:distribution}]
For brevity, let $\N_{x,\alpha}(\tau)$ be the proportion of $d\in\D_\alpha(x)$ for which $h(d)\geq \frac{2e^\gamma \sqrt x}{\log x}\tau$. Then, for any $k>1$, we have that 
\begin{align*}
k\int_0^\infty t^{k-1}\N_{x,\alpha}(t)\diff t &=\frac{1}{|\D_\alpha(x)|}\sum_{d\in\D_\alpha(x)}\int_0^{e^{-\gamma}h(d)\log x/(2\sqrt x)}kt^{k-1}\diff t
\\&=\left(\frac{2e^\gamma \sqrt x}{\log x}\right)^{-k}\frac{1}{|\D_\alpha(x)|}\sum_{d\in\D_\alpha(x)}h(d)^k.
\end{align*}
Therefore, for large $k\leq \frac{(1-2\alpha)^2}{75}\log x/(\log_2 x)^2$, combining the fact that $I_{k,1}(\alpha)\sim 2^k/(4k^2)$, our \nameref{coro_equiv}, and \eqref{eq:card_D}, one finds that
\[\int_0^\infty t^{k-1}\N_{x,\alpha}(t)\diff t=(\log k)^k\exp\left(\frac{k}{\log k}\left(C_0-1+\O\left(\frac{1}{\log k}\right)\right)\right).\]
Thus, the proof of Theorem 1.5 of \cite{lamzouri_large_2017} shows that uniformly in the range $B\leq \tau \leq \log_2 x-2\log_3 x+2\log(1-2\alpha)-\log(75)$ for some constant $B>0$:
\[\N_{x,\alpha}(\tau)=\exp\left(-\frac{e^{\tau-C_0}}\tau\left(1+\O\left(\frac{1}{\sqrt\tau}\right)\right)\right).\]
\end{proo}
\section{The number of Hooley's quadratic fields with a bounded class number}\label{sec:number_of}
Recall that $\F_\alpha(h)$ is the number of discriminants in $\D_\alpha$ with class number $h$. The first step to study $\sum_{h\leq H}\F_\alpha(h)$, for a large integer $H$, is to show that we can restrict our attention to "small" discriminants. More precisely, we will show that the main contribution to this sum comes from discriminants $d\leq X:=H^2(\log H)^2(\log\log H)^{4}$.

Thanks to Tatuzawa refinement of Siegel's inequality \cite{tatuzawa}, we know that $L(1,\chi_d)\geq 1/(\log d)$ with at most one exception. Therefore for $d\in\D_\alpha$, the class number formula implies that $h(d)\geq \sqrt d/(\log d)^2$, with at most one exception. Thus for $d\in\D_\alpha$, we have that if $h(d)\leq H$ then $d\leq H^3$, with at most one exception. This yields
\[\numberthis\label{eq:intermediaire}\sum_{h\leq H}\F_\alpha(h)=\sum_{\substack{d\in\D_\alpha(H^3)\\ h(d)\leq H}}1+\O(1).\]

Now that we only have to study discriminants $\leq H^3$, we want to refine our lower bound on $L(1,\chi_d)$, discarding a negligible set of discriminants. We will use Granville and Soundararajan's \cite{gville} result on $L(1,\chi)$, which they proved using zero-density estimates together with the large sieve.
\begin{prop}[Proposition 2.2 of \cite{gville}]\label{formule:L}
Let $A>2$ be fixed. Then, for all but at most $Q^{2/A+o(1)}$ primitive characters $\chi\mod q$ with $q\leq Q$ we have
\[L(1,\chi)=\prod_{p\leq (\log Q)^A}\left(1-\frac{\chi(p)}p\right)^{-1}\left(1+\O\left(\frac{1}{\log\log Q}\right)\right).\]
\end{prop}
We choose $A=7$, so that \nameref{formule:L} yields that for all but at most $H$ discriminants $ d\leq H^3$, we have 
\[L(1,\chi_d)=\prod_{p\leq (3\log H)^7}\left(1-\frac{\chi_d(p)}{p}\right)^{-1}(1+o(1))\geq \frac{C}{\log\log H},\]
for some $C>0$. The Class Number Formula implies that for all but at most $H$ discriminants $d\in\D_\alpha(H^3)$, $h(d)\geq C\sqrt d/(\log_2 H\log d)$. Therefore for all but at most $H$ discriminants $d\in\D_\alpha(H^3)$, $h(d)\leq H$ implies that $d\leq X$. This yields, putting this back in \eqref{eq:intermediaire}:
\[\numberthis\label{eq:tatuzawa}\sum_{h\leq H}\F_\alpha(h)=\sum_{\substack{d\in\D_\alpha(X)\\ h(d)\leq H}}1+\O(H).\]

We are now ready to prove our theorem, following Dahl-Lamzouri's \cite{dahl_distribution_2016} and Soundararajan's \cite{soundararajan} papers.
\begin{proo}[of \nameref{theo:number_F}]
We introduce the smoothing function
\[I_{c,\lambda,N}(y):=\frac{1}{2\pi i}\int_{c-i\infty}^{c+i\infty} y^s\left(\frac{e^{\lambda s}-1}{\lambda s}\right)^N\frac{\diff s}s,\]
where $c=1/\log H$, $N$ is a positive integer, and $0<\lambda\leq 1$ is a real number. By Perron's formula, we know that (see $(4.19)$ of \cite{dahl_distribution_2016})
\[\numberthis\label{formule_smooth}I_{c,\lambda,N}(y)=\left\{\begin{array}{cl}
=1&\text{if }y\geq 1,
\\\in[0,1]&\text{if }e^{-\lambda N}\leq y<1,
\\=0&\text{if }0<y<e^{-\lambda N}.
\end{array}\right.\]
According to \eqref{eq:tatuzawa} and \eqref{formule_smooth}, we have that 
\[\numberthis\label{eq:main_Fh}\sum_{h\leq H}\F_\alpha(h)\leq \frac{1}{2\pi i}\int_{c-i\infty}^{c+i\infty}\sum_{d\in\D_\alpha(X)}\frac{H^s}{h(d)^s}\left(\frac{e^{\lambda s}-1}{\lambda s}\right)^N\frac{\diff s}s+\O(H)\leq \sum_{h\leq e^{\lambda N}H}\F_\alpha(h).\]

Thanks to \nameref{theo-hd} we know that uniformly in $|s|\leq T:=(\log X)^{1/3-\varepsilon}$, $\varepsilon>0$ small,
\[\sum_{d\in\D_\alpha(X)}h(d)^{-s}=\mathfrak I_1+\mathfrak I_2+\O(\sqrt X(\log\log X)^3T),\]
where
\[\mathfrak I_1=\mathfrak I_1(s):=\frac{1}{2\zeta(2)}\int_2^X\frac{t^{\frac{-s-1}2}}{\log(t)^{-s}}\EE(L(1,\XX)^{-s})I_{-s,1}(\alpha)\log^2(t)\diff t,\]
and
\[\mathfrak I_2=\mathfrak I_2(s):=\frac{1}{2\zeta(2)}\int_2^X\frac{t^{\frac{-s-1}2}}{\log(t)^{-s}}\EE(L(1,\XX)^{-s})\phi(-s)I_{-s,0}(\alpha)\log(t)\diff t.\]

Since $h(d)\geq 1$, $|e^{\lambda s}-1|\leq 3$ for large $H$ and $|\D_\alpha(X)|\ll \sqrt X\log^2 X$ by \eqref{eq:card_D}, the contribution of $s$ such that $\Re(s)=c$ and $|s|>T$ in the integral of \eqref{eq:main_Fh} is 
\[\ll|\D_\alpha(X)|\left(\frac{3}{\lambda}\right)^N\int_{\substack{|s|>T\\\Re(s)=c}}\frac{|\diff s|}{|s|^{N+1}}\ll \frac{X^{1/2}\log^2X}N\left(\frac{3}{\lambda T}\right)^N.\]
For large $H$ we also have that $|(e^{\lambda s}-1)/(\lambda s)|\leq 4$, which is easily seen by separating the cases $|\lambda s|>1$ and $|\lambda s|\leq 1$. Thus, the integral term of \eqref{eq:main_Fh} is equal to 
\begin{align*}
\numberthis\label{eq:Fh_int}\frac{1}{2\pi i}\int_{\substack{\Re(s)=c\\|s|\leq T}}H^s\left(\mathfrak I_1+\mathfrak I_2\right)\left(\frac{e^{\lambda s}-1}{\lambda s}\right)^N\frac{\diff s}s+R,
\end{align*}
where 
\begin{align*}
R&\ll \frac{X^{1/2}\log^2X}N\left(\frac{3}{\lambda T}\right)^N+\frac{4^NT^2}cX^{1/2}(\log\log X)^3.
\end{align*}
Let $\lambda=3e^{1/\varepsilon}/T$ and $N=\left\lfloor\varepsilon\log_2H\right\rfloor$, so that 
\[R\ll_\varepsilon H(\log H)^{8/3}.\]

Note that on the line $\Re(s)=c$, $I_{-s,1}(\alpha),I_{-s,0}(\alpha),\EE(L(1,\XX)^{-s}),\EE(L(1,\XX)^{-s})\phi(-s)\ll 1$ (by \nameref{rema:phi}). Therefore we may extend the segment of integration to the whole line $\Re(s)=c$ in \eqref{eq:Fh_int}, the introduced error being 
\begin{align*}
\ll\left(\frac{3}\lambda\right)^N&\int_{\substack{|s|>T\\ \Re(s)=c}}\int_2^Xt^{\frac{-1}2}\log^2 t\diff t\frac{|\diff s|}{|s|^{N+1}}\ll \frac{X^{1/2}\log^2 X}{N}\left(\frac{3}{\lambda T}\right)^N\ll_\varepsilon H(\log H)^{8/3}.
\end{align*}

Thus, the middle term of \eqref{eq:main_Fh} is equal to
\begin{align*}
\frac{1}{2\pi i}\int_{\substack{\Re(s)=c}}H^s(\mathfrak I_1+\frak I_2)\left(\frac{e^{\lambda s}-1}{\lambda s}\right)^N\frac{\diff s}s+\O_\varepsilon(H(\log H)^{8/3}).
\end{align*}
Now we split the integral in two, and study each term separately. By switching the integrals, one shows that 
\[\numberthis\label{eq:doubleint}\frac{1}{2\pi i}\int_{\substack{\Re(s)=c}}H^s\mathfrak I_1\left(\frac{e^{\lambda s}-1}{\lambda s}\right)^N\frac{\diff s}s=\frac{1}{\zeta(2)}\int_{1/2}^{1/2+\alpha}\EE\left[\int_2^X\frac{\log^2 t}{2\sqrt t}I_{c,\lambda,N}\left(\frac{Hu\log t}{L(1,\XX)\sqrt t}\right)\diff  t\right]\left(u-\frac12\right)\diff u.\]

Put $Y:=H/L(1,\XX)$ and note that $f:t\mapsto \frac{\sqrt t}{\log t}$ is strictly increasing (thus invertible) on $(10,\infty)$: we denote its inverse by $g$. With these notations, if $1/2\leq u\leq \alpha+1/2$, we may use \eqref{formule_smooth} to find that
\[I_{c,\lambda,N}\left(\frac{Hu\log t}{\sqrt tL(1,\XX)}\right)=\left\{\begin{array}{cl}
=1&\text{if }uY\geq f(t),
\\=0&\text{if }uYe^{\lambda N}<f(t),
\\\in[0,1]&\text{otherwise}.
\end{array}\right.\]
Thus, if we let $\mathfrak F(t)=\sqrt t\left((\log t)^2+\O(\log t)\right)$ be a primitive of $t\mapsto(\log t)^2/(2\sqrt t)$, we get that
\begin{align*}
\numberthis\label{eq:somme1}\int_2^X\frac{\log^2 t}{2\sqrt t}I_{c,\lambda,N}\left(\frac{Hu\log t}{\sqrt tL(1,\XX)}\right)\diff t&=\min(\mathfrak F(X), \mathfrak F(g(uY)))
\\&+\O\left(\mathfrak F(g(uYe^{\lambda N}))-\mathfrak F(g(uY))+1\right).
\end{align*}
Note that $g(t)=4t^2(\log t+\O(\log \log t))^2$, for $t$ large enough, so that $\mathfrak F(g(t))=8t\log^3(t)+\O(t\log^2(t)\log_2(t))$. Furthermore for $1/2\leq u\leq 1$, if $g(uY)>X$, then $uY>\sqrt X/\log X$ and hence $L(1,\XX)\ll (\log H)^{-2}$. By \nameref{theo:phi_esti}, this happens with probability $\ll \exp(-H)$. Thus the Cauchy-Schwarz inequality yields
\begin{align*}
\EE\bigr[\min(\mathfrak F(g(uY)),\mathfrak F(X))\bigl]&=\EE(\mathfrak F(g(uY)))+\EE[(\frak F(X)-\frak F(g(uY)))\mathfrak 1_{X<g(uY)}]
\\&=\EE(\mathfrak F(g(uY)))+\O(\EE[\frak F(g(uY))\mathfrak 1_{X<g(uY)}])
\\&=\EE(\mathfrak F(g(uY)))+\O(\sqrt{\EE[\frak F(g(uY))^2]\PP(X<g(uY))})
\\&=8u\EE(L(1,\XX)^{-1})H\log^3H+\O(1).
\end{align*}
Moreover
\begin{align*}
\EE\bigg[\mathfrak F(g(uYe^{\lambda N}))-\mathfrak F(g(uY))\bigg]&\ll(e^{\lambda N}-1)H\log^3 H+H\log^2H\log_2H
\\&\ll_\varepsilon H(\log H)^{8/3+2\varepsilon}.
\end{align*}
Putting these back in \eqref{eq:somme1} and \eqref{eq:doubleint}, we get that 
\begin{align*}
&\frac{1}{2\pi i}\int_{\substack{\Re(s)=c}}H^s\mathfrak I_1\left(\frac{e^{\lambda s}-1}{\lambda s}\right)^N\frac{\diff s}s
\\&\numberthis\label{eq_I1_firstesti}\qquad=\frac{8}{\zeta(2)}\EE(L(1,\XX)^{-1})\int_{1/2}^{\alpha+1/2} u\left(u-\frac12\right)\diff u \times H\log^3 H+\O_\varepsilon(H(\log H)^{8/3+2\varepsilon}).
\end{align*}

Now we need to show that the contribution of $\frak I_2$ goes into the error term. By definition, 
\[\EE(L(1,\XX)^{-s})\phi(-s)=\sum_{p}\prod_{q\neq p}\EE\left(1-\frac{\XX(p)}p\right)^s\left[c_1(p)\left(1-\frac1p\right)^{s}+c_2(p)\left(1+\frac1p\right)^s+c_3(p)\right],\]
for some coefficients $c_i(p)\ll \log p/p^2$. We introduce the random variables $\YY_p(q)$, $p,q$ primes, defined by the following law: $\YY_p(q)$ and $\XX(q)$ have the same law if $p\neq q$ and $\PP(\YY_p(p)=0)=1$. We then put, for $n\geq 2$: $\YY_p(n):=\prod_{q^e||n}\YY_p(q)^e$ and $\YY_p(1)=1$. With this formalism, we have 
\[\prod_{q\neq p}\EE\left(1-\frac{\XX(p)}p\right)^s=\EE(L(1,\YY_p)^{-s}).\]
Then,
\begin{align*}
\frac{1}{2\pi i}\int_{\substack{\Re(s)=c}}&H^s\mathfrak I_2\left(\frac{e^{\lambda s}-1}{\lambda s}\right)^N\frac{\diff s}s=
\\\frac{1}{\zeta(2)}\sum_p\Biggl\{&c_1(p)\int_{1/2}^{1/2+\alpha}\EE\left[\int_2^X\frac{\log t}{2\sqrt t}I_{c,\lambda,N}\left(\frac{Hu\log t}{L(1,\YY_p)\sqrt t}\left(1-\frac1p\right)\right)\diff  t\right]\diff u
\\+&c_2(p)\int_{1/2}^{1/2+\alpha}\EE\left[\int_2^X\frac{\log t}{2\sqrt t}I_{c,\lambda,N}\left(\frac{Hu\log t}{L(1,\YY_p)\sqrt t}\left(1+\frac1p\right)\right)\diff  t\right]\diff u
\\+&c_3(p)\int_{1/2}^{1/2+\alpha}\EE\left[\int_2^X\frac{\log t}{2\sqrt t}I_{c,\lambda,N}\left(\frac{Hu\log t}{L(1,\YY_p)\sqrt t}\right)\diff  t\right]\diff u\Biggr\}.
\end{align*}
We now deal with each of these integrals in a similar way to that of the contribution of $\frak I_1$. We only have to consider $\tilde{\frak F}(t)=\sqrt t(\log t+\O(1))$, a primitive of $t\mapsto \log t/(2\sqrt t)$, instead of $\frak F$. As a result, we find that each of the three integrals in $u$ between curly brackets above is $\ll H\log^2H$, and hence 
\[\frac{1}{2\pi i}\int_{\substack{\Re(s)=c}}H^s\mathfrak I_2\left(\frac{e^{\lambda s}-1}{\lambda s}\right)^N\frac{\diff s}s\ll \sum_p (c_1(p)+c_2(p)+c_3(p)) H\log^2H\ll H\log^2 H\]
since $c_i(p)\ll \log p/p^2$.

All in one, putting this with \eqref{eq_I1_firstesti} in \eqref{eq:Fh_int}, we find that the middle term of \eqref{eq:main_Fh} is equal to
\[\frac{8}{\zeta(2)}\EE(L(1,\XX)^{-1})\int_{1/2}^{\alpha+1/2} u\left(u-\frac12\right)\diff u \times H\log^3 H+\O_\varepsilon(H(\log H)^{8/3+2\varepsilon}).\]
Since $\EE(L(1,\XX)^{-1})=\prod_p\left(1+\frac{p-1}{p^2(p+1)}\right)$ the previous expression becomes
\[2\frac{\alpha^2(4\alpha+3)}{3}\prod_p\left(1-\frac{2p-1}{p^4}\right)H\log^3 H+\O_\varepsilon(H(\log H)^{8/3+2\varepsilon}).\]
This concludes the proof, changing $e^{\lambda N}H$ by $H$ in the right hand side of \eqref{eq:main_Fh}.
\end{proo}
\bibliographystyle{abbrv}
\bibliography{biblio}
\labo
\end{document}